\documentclass[10pt]{article}
\usepackage{latexsym,amsmath,amssymb,graphics,amscd}
\usepackage{multirow}
\usepackage{soul}
\usepackage{booktabs}
\usepackage{tabularx}
\usepackage[hang,footnotesize]{caption}
\usepackage[pdftex]{graphicx}
\usepackage{graphicx}
\usepackage{overpic}
\usepackage{appendix}
\usepackage{cases}
\usepackage{color}
\usepackage{tikz}
\usetikzlibrary{decorations.pathreplacing}
\addtocounter{footnote}{1}
\makeatletter
\@addtoreset{table}{bsection}
\def\thetable{\thesection.\@arabic\c@table}
\def\fps@table{h, t}
\makeatother

\newtheorem{theorem}{Theorem}[section]
\newtheorem{definition}[theorem]{Definition}
\newtheorem{lemma}[theorem]{Lemma}

\newtheorem{proposition}[theorem]{Proposition}

\newtheorem{corollary}[theorem]{Corollary}

\newsavebox{\savepar}

\makeatletter

\begin{document}

\title{\textbf{Optimal nonlinear information processing capacity in delay-based reservoir computers}}
\author{Lyudmila Grigoryeva$^{1}$, Julie Henriques$^{1, 2}$, Laurent Larger$^{3}$, and Juan-Pablo Ortega$^{4, \ast}$}
\date{}
\maketitle

\begin{abstract}
Reservoir computing  is a recently introduced brain-inspired machine learning paradigm capable of excellent performances in the processing of empirical data. We focus in a particular kind of time-delay based reservoir computers that have been physically implemented using optical and electronic systems and have shown unprecedented data processing rates. Reservoir computing is well-known for the ease of the associated training scheme but also for the problematic sensitivity of its performance to architecture parameters. This article addresses the reservoir design problem, which remains the biggest challenge in the applicability of this information processing scheme. More specifically, we use the information available regarding the optimal reservoir working regimes to construct a functional link between the reservoir parameters and its performance. This function is used to explore various properties of the device and to choose the optimal reservoir architecture, thus replacing the tedious and time consuming parameter scannings used so far in the literature.
\end{abstract}

\bigskip

\textbf{Key Words:} Reservoir computing, echo state networks, neural computing, time-delay reservoir, memory capacity, architecture optimization.

\makeatletter
\addtocounter{footnote}{1} \footnotetext{%
Laboratoire de Math\'{e}matiques de Besan\c{c}on, UMR CNRS 6623, Universit\'{e} de Franche-Comt\'{e}, UFR des
Sciences et Techniques. 16, route de Gray. F-25030 Besan\c{c}on cedex. France. {\texttt{Lyudmyla.Grygoryeva@univ-fcomte.fr} }}
\makeatother
\makeatletter
\addtocounter{footnote}{1} \footnotetext{%
Cegos Deployment. 11, rue Denis Papin. F-25000 Besan\c{c}on. {\texttt{jhenriques@deployment.org} }}
\makeatother
\makeatletter
\addtocounter{footnote}{1} \footnotetext{%
FEMTO-ST, UMR CNRS  6174, Optics Department, Universit\'{e} de Franche-Comt\'{e}, UFR des
Sciences et Techniques. 15, Avenue des Montboucons. F-25000 Besan\c{c}on cedex. France. {\texttt{Laurent.Larger@univ-fcomte.fr} }}
\makeatother

\makeatletter
\addtocounter{footnote}{1} \footnotetext{%
Corresponding author. Centre National de la Recherche Scientifique, Laboratoire de Math\'{e}matiques de Besan\c{c}on, UMR CNRS 6623, Universit\'{e} de Franche-Comt\'{e}, UFR des
Sciences et Techniques. 16, route de Gray. F-25030 Besan\c{c}on cedex.
France. {\texttt{Juan-Pablo.Ortega@univ-fcomte.fr} }}
\makeatother

\medskip

\medskip

\medskip

The increase in need for information processing capacity, as well as the physical limitations of the Turing or von Neumann machine methods implemented in most computational systems, has motivated  the search for new brain-inspired solutions some of which present an outstanding potential. An important direction in this undertaking is based on the use of the intrinsic information processing abilities of dynamical systems~\cite{Crutchfield2010} which opens the door to high performance physical realizations whose behavior is ruled by these structures~\cite{Caulfield2010, Woods2012}.

The contributions in this paper take place in a specific implementation of this idea that is obtained as a melange of a recently introduced machine learning paradigm known under the name of {\bf reservoir computing (RC)}~\cite{jaeger2001, Jaeger04, maass1, maass2, Crook2007, verstraeten, lukosevicius} with a realization based on the sampling of the solution of a time-delay differential equation~\cite{Rodan2011, gutierrez2012}.  We  refer to this combination as {\bf time-delay reservoirs (TDRs)}. Physical implementations of this scheme carried out with dedicated hardware are already available and have shown excellent performances in the processing of empirical data:  spoken digit recognition~\cite{jaeger2, Appeltant2011, Larger2012, Paquot2012, photonicReservoir2013}, the NARMA model identification task~\cite{Atiya2000, Rodan2011}, continuation of chaotic time series, and volatility forecasting~\cite{GHLO2012}. A recent example that shows the potential of this combination are the results in~\cite{photonicReservoir2013} where an optoelectronic implementation of a TDR is capable of achieving the lowest documented error in the speech recognition task at unprecedented speed in an experiment design in which digit and speaker recognition are carried out in parallel.

A major advantage of RC is the linearity of its training scheme. This choice makes its implementation easy when compared to more traditional machine learning approaches like recursive neural networks, which usually require the solution of convoluted and sometimes ill-defined optimization problems. In exchange, as it can be seen in most of the references quoted above, the system performance is not robust with respect to the choice of the parameter values $\boldsymbol{\theta} $ of the nonlinear kernel used to construct the RC (see below). More specifically, small deviations from the optimal parameter values can seriously degrade the performance and moreover, the optimal parameters are highly dependent on the task at hand. This observation makes the kernel parameter optimization a very important step in the RC design and has motivated the introduction of alternative parallel-based architectures~\cite{pesquera2012, GHLO2012} to tackle this difficulty.  

The main contribution of this paper is the introduction of an approximated model that, to our knowledge, provides the first rigorous analytical description of the delay-based RC performance. This powerful theoretical tool can be used to systematically study the delay-based RC properties and to replace the trial and error approach in the choice of architecture parameters by well structured optimization problems. This method simplifies enormously the implementation effort and sheds new light on the mechanisms that govern this information processing technique.

TDRs are based on the interaction of the time-dependent input signal $z(t) \in \mathbb{R}   $ that we are interested in with the solution space of a time-delay differential equation of the form
\begin{equation}
\label{time delay equation 2}
\dot{x} (t)= - x (t)+f(x (t- \tau), I(t), \boldsymbol{\theta}),
\end{equation}
where $f  $ is a nonlinear smooth function (we call it {\bf nonlinear kernel}) that depends on the $K$ parameters in the vector $\boldsymbol{\theta} \in \mathbb{R} ^K $, $\tau >0$ is the {\bf delay}, $x (t) \in \mathbb{R} $, and $I (t)  \in \mathbb{R}$ is obtained using a temporal multiplexing over the delay period of the input signal $z (t) $ that we explain later on. We note that, even though the differential equation takes values in the real line, its solution space is infinite dimensional since an entire function $x \in C^{1}([- \tau, 0], \mathbb{R}) $ needs to be specified in order to initialize it. The choice of nonlinear kernel is determined by the intended physical implementation of the computing system; we focus on two parametric sets of kernels that have already been explored in the literature, namely, the Mackey-Glass~\cite{mackey-glass:paper} and the Ikeda~\cite{Ikeda1979} families. These kernels were used for reservoir computing purposes in the RC electronic and optic realizations in~\cite{Appeltant2011} and~\cite{Larger2012}, respectively.   

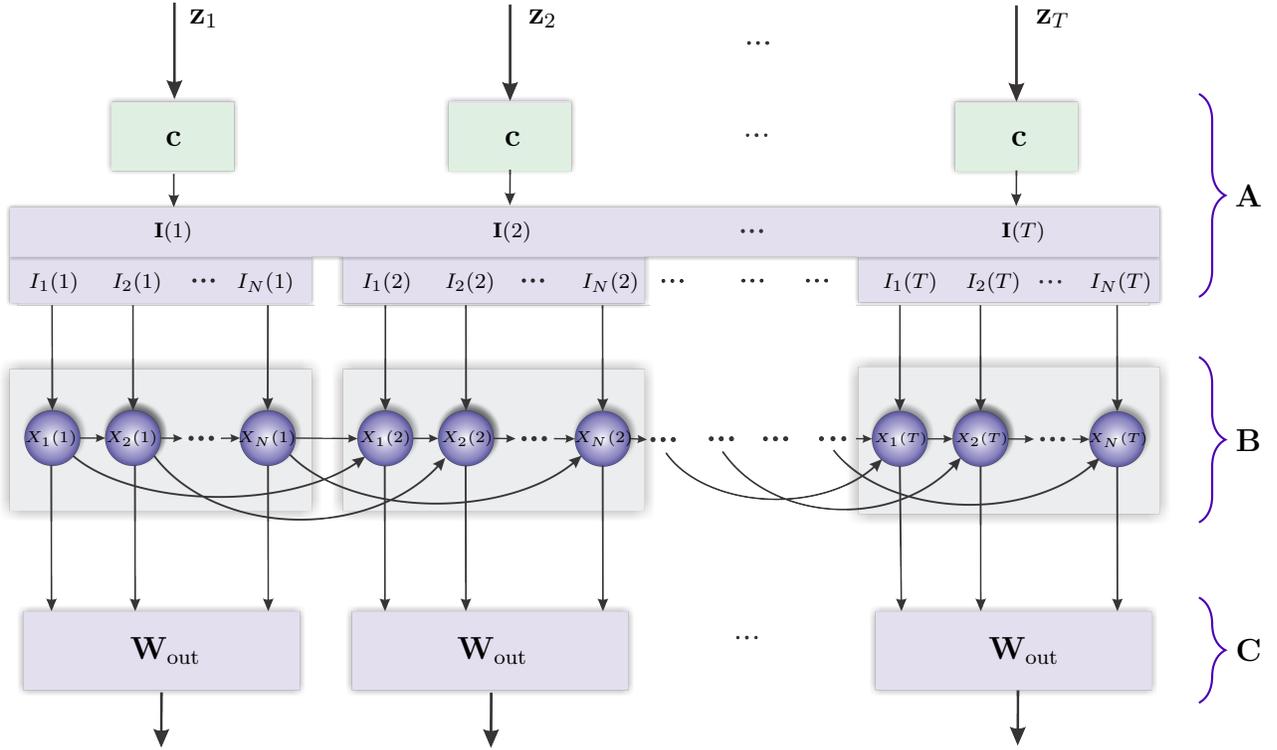
\begin{figure}[!htp]
\hspace{-0.5cm}
\begin{tikzpicture}
{\footnotesize
\begin{overpic}[scale = 0.35,tics = 5]{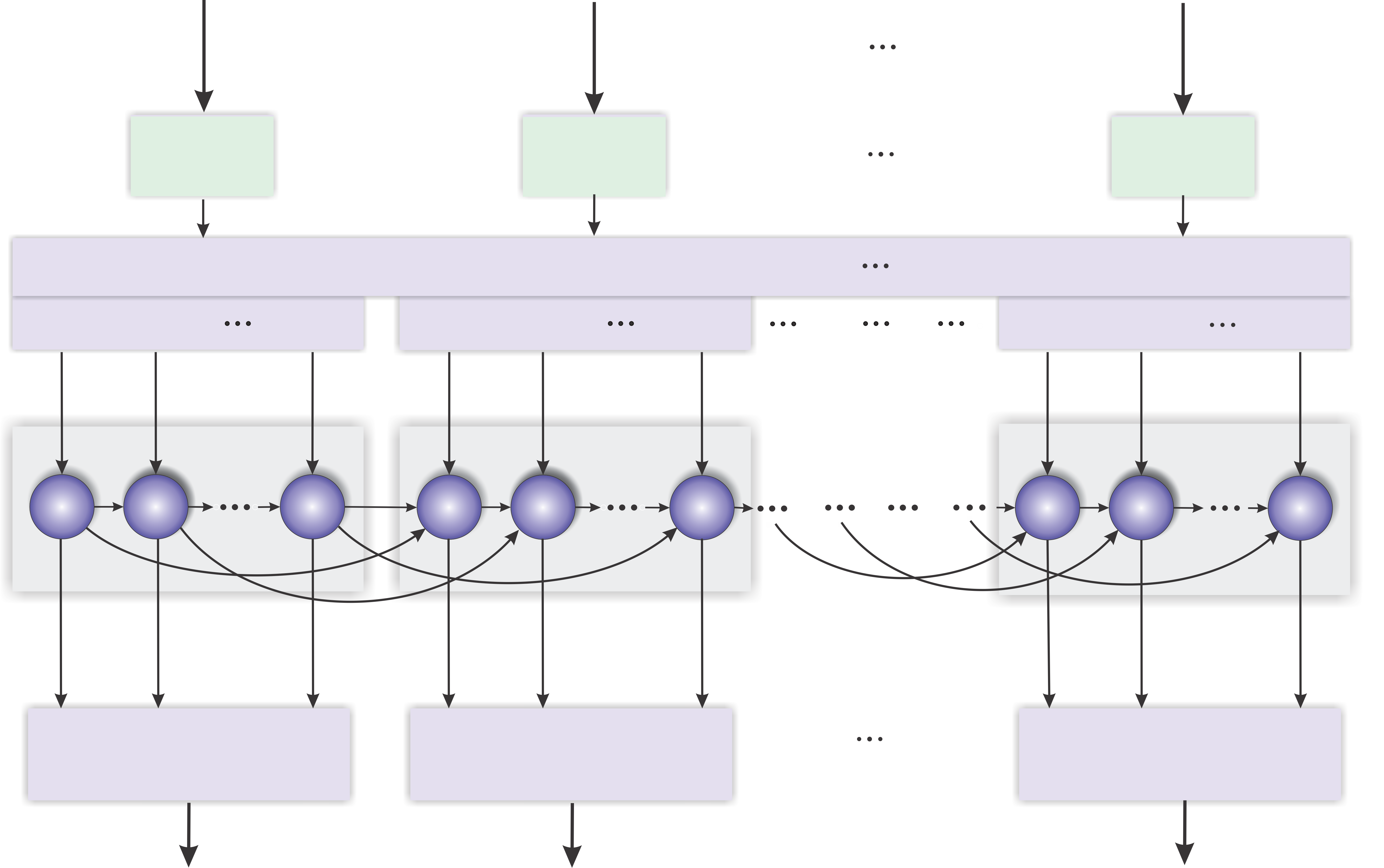}
\put (14,50.7) {\large$\bf c$}
\put (42.5,50.7) {\large${\bf c} $}
\put (85.1,50.7) {\large$ {\bf c} $}

\put (2.2,25.7) {\fontsize{6}{100}$X_{1}(1)$}
\put (9,25.7) {\fontsize{6}{100}$X_{2}(1)$}
\put (20.1,25.7) {\fontsize{6}{100}$X_{N}(1)$}

\put (30.3,25.7) {\fontsize{6}{100}$X_{1}(2)$}
\put (37.2,25.7) {\fontsize{6}{100}$X_{2}(2)$}
\put (48.4,25.7) {\fontsize{6}{100}$X_{N}(2)$}

\put (73.6,25.7) {\fontsize{5}{100}$X_{1}(T )$}
\put (80.4,25.7) {\fontsize{5}{100}$X_{2}(T )$}
\put (91.6,25.7) {\fontsize{5}{100}$X_{N}(T )$}


\put (16,61) {\large$\mathbf{z}_1$}
\put (44.5,61) {\large$\mathbf{z}_2$}
\put (87.2,61) {\large$\mathbf{z}_{T}$}

\put (13,43) {$\mathbf{I}(1)$}

\put (41.5,43) {$\mathbf{I}(2)$}

\put (84.3,43) {$\mathbf{I}(T)$}

\put (2.5,38.7) {$I_1(1)$}
\put (9.5,38.7) {$I_2(1)$}
\put (20,38.7) {$I_N(1)$}

\put (30.5,38.7) {$I_1(2)$}
\put (37.5,38.7) {$I_2(2)$}
\put (49,38.7) {$I_N(2)$}

\put (74.3,38.7) {$I_1(T)$}
\put (81.3,38.7) {$I_2(T)$}
\put (91.7,38.7) {$I_N(T)$}

\put (83.2,7.5) {\large ${\bf W}_{\rm out}$}
\put (11,7.5) {\large ${\bf W}_{\rm out}$}
\put (38.5,7.5) {\large ${\bf W}_{\rm out}$}

\end{overpic} }
\draw [decorate,decoration={brace,amplitude=10pt,mirror,raise=4pt},yshift=0pt,red!30!blue,thick,solid]
(0,0.6) -- (0,2) node [black,midway,xshift=0.8cm] {\large
{\bf C}};
\draw [decorate,decoration={brace,amplitude=10pt,mirror,raise=4pt},yshift=0pt,red!30!blue,thick,solid]
(0,3) -- (0,5.2) node [black,midway,xshift=0.8cm] {\large
{\bf B}};
\draw [decorate,decoration={brace,amplitude=10pt,mirror,raise=4pt},yshift=0pt,red!30!blue,thick,solid]
(0,6) -- (0,8.7) node [black,midway,xshift=0.8cm] {\large
{\bf A}};
\end{tikzpicture}\vspace{0.7cm}
\caption{Neural diagram representing the architecture of the time-delay reservoir (TDR) and the three modules of the reservoir computer (RC): A is the input layer, B is the time-delay reservoir, and C is the readout layer.}
\label{RC figure}
\end{figure}

In order to visualize the TDR construction using a neural networks approach it is convenient, as in~\cite{Appeltant2011, gutierrez2012}, to consider the Euler time-discretization of~\eqref{time delay equation 2} with integration step $d:= \tau/N $, namely,
\begin{equation}
\label{euler discretization}
\frac{x (t)- x (t- d)}{d}=- x (t)+ f(x (t- \tau), I(t), \boldsymbol{\theta}).
\end{equation}
The design starts with the choice of a number $N\in \mathbb{N}$ of {\bf  virtual neurons} and of an adapted {\bf input mask} ${\bf c}\in  \mathbb{R}^N $. 
Next, the input signal $z (t) $ at a given time $t$ is multiplexed over the delay period by setting $\mathbf{I}(t):={\bf c} z (t) \in \mathbb{R}^N$ (see module A in Figure~\ref{RC figure}). We  then organize it, as well as the solutions of~\eqref{euler discretization}, in {\bf neuron layers} $\mathbf{x} (t)  $ parametrized by a discretized  time $t\in \mathbb{Z} $ by setting
\begin{equation*}
x_i (t) := x (t \tau- (N-i)d),  \quad I_i (t) := I (t \tau- (N-i)d),\quad i \in \{1, \ldots, N\}, \quad t\in  \mathbb{Z},
\end{equation*}
where $x_i (t) $ and $I_i (t) $ stand for the $i $th-components of the vectors $\mathbf{x} (t)  $ and $\mathbf{I} (t)  $, respectively, with $t \in  \mathbb{Z} $. We  say  that $x _i (t) $ is the {\bf  $i$th neuron value of the $t $th layer of the reservoir} and $d $ is referred to as the {\bf  separation between neurons}. With this convention, the solutions of~\eqref{euler discretization} are described by the following recursive relation:
\begin{equation}
\label{recursion euler}
x _i(t):= e^{-\xi}x_{i-1}(t)+(1-e^{-\xi})f(x _i(t-1), I _i(t), \boldsymbol{\theta}),\enspace \mbox{with} \enspace x _0(t):= x _N(t-1), \enspace \mbox{and} \enspace \xi:=\log (1+ d),
\end{equation}
that shows how, as depicted in module B in Figure~\ref{RC figure}, any  neuron value is a convex linear combination of the previous neuron value in the same layer and a nonlinear function of both the same neuron value in the previous layer and the input. The weights of this combination are determined by the separation between neurons; when the distance $d$ is small, the neuron value $x_i(t)$  is mainly influenced by the previous neuron value $x_{i-1}(t)$, while large distances between neurons give predominance to the previous layer and foster the input gain. The recursions~\eqref{recursion euler} uniquely determine a smooth map
$F: \mathbb{R}^N\times \mathbb{R}^N \times \mathbb{R} ^K \rightarrow \mathbb{R}^N $ 
that specifies the neuron values as a recursion on the neuron layers via an expression of the form 
\begin{equation}
\label{vector discretized reservoir main}
\mathbf{x}(t) = F( \mathbf{x} (t-1), \mathbf{I} (t), \boldsymbol{ \theta }) , 
\end{equation}
where $F$ is constructed out of the nonlinear kernel map $f$ that depends on the $K$ parameters in the vector $\boldsymbol{\theta}$; $F$ is referred to as the {\bf reservoir map}.

The construction of the TDR computer is finalized by connecting, as in Module C of Figure~\ref{RC figure}, the reservoir output to a linear readout ${\bf W}_{{\rm out}} \in \mathbb{R} ^N $ that is calibrated using a training sample by minimizing the associated task mean square error via a linear regression. We will refer to the module B in Figure~\ref{RC figure} as the {\bf reservoir} or the {\bf time-delay reservoir} (TDR) and to the collection of the three modules as the {\bf reservoir computer} (RC) or the TDR computer. A TDR based on the direct sampling of the solutions of~\eqref{time delay equation 2} will be called a {\bf continuous time TDR} and those based on the recursion~\eqref{vector discretized reservoir main} will be referred to as {\bf discrete time TDRs}.

\medskip


\vspace{1cm}

As we already mentioned, the performance of the RC for a given task is much dependent on the value of the kernel parameters $\boldsymbol{\theta} $ and, in some cases, on the entries of the input mask ${\bf c}  $ used for signal multiplexing. The optimal parameters $\boldsymbol{\theta} $ are usually determined by trial and error or using computationally costly systematic scannings that are by far the biggest burden at the time of adapting the RC to a new task. In this paper we construct an approximate model that we use to allows us to establish a functional link between the RC performance and the parameters $\boldsymbol{\theta} $ and the input mask values ${\bf c} $. Given a specific task, this explicit expression can be used to find appropriate parameter and mask values by solving a well structured and algorithmically convenient optimization problem that readily provides them.

The construction of this approximated formula is based on the observation that the optimal RC performance  is always obtained when the TDR is working in a {\bf stable unimodal regime}, that is, the reservoir is initialized at a stable equilibrium of the autonomous system ($I  (t)=0$) associated to~\eqref{time delay equation 2} and the mean and variance of the input signal $I (t)  $ are designed using the input mask ${\bf c}  $ so that the reservoir output remains around it and does not visit other stable equilibria or  dynamical elements. In the next section we provide empirical and theoretical arguments for this claim. The performance measures that we consider in our study are the nonlinear memory capacities introduced in~\cite{dambre2012} as a generalization of the linear concept proposed in~\cite{Jaeger:2002, White2004, Ganguli2008, Hermans2010}.

\section{Results}
\label{Results} 

\subsection{Optimal performance: stability and unimodality}
\label{Optimal performance: stability and unimodality}

\paragraph{Stability and the reservoir defining properties.} 

The estimations of the RC performance using the nonlinear memory capacity that we present later on, consist of approximating the reservoir by its partial linearization at the level of the delayed self feedback term and of respecting the nonlinearity in the input injection. This approach is only acceptable when the optimal dynamical regime that we are interested in, remains close to a given point. A natural candidate for such qualitative behavior could be obtained by initializing the reservoir at an asymptotically stable equilibrium of the autonomous system associated to~\eqref{time delay equation 2} and by controlling the mean and the variance of the input signal $I (t)  $ so that the reservoir output remains close to it. 

There is both theoretical and empirical evidence that suggests that optimal performance is obtained when working in a statistically stationary regime around a stable equilibrium. Indeed, one of the defining features of RC, namely the {\bf echo state property} is materialized for general RCs by enforcing that the spectral radius of the internal connectivity matrix of the reservoir is smaller than one~\cite{maass1,Jaeger04,lukosevicius}, which is the critical stability value for a quiescent state of the network when operating  autonomously  (without external injected information). It is well-known  that the translation of this condition for TDRs implies parameter settings that ensure the existence of a stable state of~\eqref{time delay equation 2} when $I (t) $ is set to zero. This feature typically relates to gains of the feedback smaller than the Hopf threshold of the delay dynamics or, equivalently, to a sufficiently low feedback rate so that self-sustained oscillations are avoided. 

Asymptotic stability is closely related with the so-called {\bf fading memory property} \cite{Boyd1985,maass1}: the impact of any past injected input necessarily vanishes after a
transient whose duration is typically of the order of the absolute value of the inverse of the smallest negative real part in the Lyapunov exponents. When the feedback gain is set too close to zero, the RC does not exhibit a long enough transient and thus presents an intrinsic memory that is too short to secure the  self mixing of the temporal information necessary for its processing. On the other hand, if the feedback gain is set too close to the instability threshold, the input information flow requires too much time to vanish and hence the fading memory property is poorly satisfied. We recall the well known fact (see Section 8.2 in~\cite{Boyd1985}) that the fading memory property can be realized by input-output systems generated by time-delay differential equations only when these exhibit a unique stable equilibrium. 

In the context of recent successful physical realizations of RC, experimental parameters are systematically chosen so that the conditions described above are satisfied. Indeed, in \cite{Appeltant2011,Larger2012} these conditions are ensured via a proper tuning of the
gain of the delayed feedback function. This approach differs from the one in~\cite{photonicReservoir2013}, where the conditions are met by choosing a laser injection current strictly smaller  but close to the lasing threshold, as well as by using a moderate feedback, which prevents eventual self sustained external cavity mode oscillations. An additional important observation suggested by all these experimental setups is the need for a nonlinearity at the level of the input injection. In \cite{Appeltant2011,Larger2012}
this feature is obtained using a strong enough input signal amplitude and via the transformation associated  to the nonlinear delayed feedback. In \cite{photonicReservoir2013}  the delayed feedback is linear but an external Mach-Zehnder modulator is used that implicitly
provides a nonlinear transformation of the input signal as it is optically seeded through the nonlinear electro-optic modulation transfer function of the Mach-Zehnder. 

\paragraph{Stability analysis of the time-delay reservoir.} Due to the central role played by stability in our discussion, we now carefully analyze various sufficient conditions that ensure that the RC is functioning in a stable regime. All the statements that follow are carefully proved in the Supplementary Material section. Consider first an equilibrium $x _0 \in \mathbb{R} $ of the continuous time model~\eqref{time delay equation 2} working in autonomous regime, that is, we set $I(t)=0 $. It can be shown using a Lyapunov-Krasovskiy-type analysis~\cite{krasovskiy:book, she:book} that the asymptotic stability of $x _0 $ is guaranteed whenever there exists an $\varepsilon >0 $  and a constant $|k_\varepsilon|<1$ such that  either
\begin{description}
\item[(i)] $f(x + x _0,0, \boldsymbol{ \theta }  )\le k_\varepsilon x + x _0 $ for all $x \in \left( - \varepsilon , \varepsilon\right) $, or
\item[(ii)] $\dfrac{f(x + x _0 ,0, \boldsymbol{ \theta }) - x_0 }{x}\le k_\varepsilon$ for all $x \in \left( - \varepsilon , \varepsilon\right) $.
\end{description} 
The first condition can be used to prove the stability of equilibria exhibited by TDRs created using concave (but not necessarily differentiable) nonlinear kernels. As to the second one, it shows that if $f$ is differentiable at $x _0 $ then this point is stable as long as  $|\partial _{x}f(x _0, 0, \boldsymbol{\theta}  )|<1 $, with $\partial _{x}f(x _0, 0, \boldsymbol{\theta}  )$ the first derivative of the nonlinear kernel $f$ in \eqref{time delay equation 2} with respect to the first argument at the point $(x _0, 0 , \boldsymbol{\theta}  )$. 

The stability study can also be carried out by working with the discrete-time approximation~\eqref{vector discretized reservoir main} of the TDR which is determined by the reservoir map  $F: \mathbb{R}^N\times \mathbb{R}^N \times \mathbb{R} ^K \rightarrow \mathbb{R}^N $. More specifically, it can be shown that $x _0 \in \mathbb{R} $ is an equilibrium of~\eqref{time delay equation 2} if and only if $\mathbf{x} _0:=(x _0, \ldots, x _0)^\top \in  \mathbb{R}^N$ is a fixed point of  \eqref{vector discretized reservoir main}. The asymptotic stability of this fixed point is ensured whenever the linearization $D_{\mathbf{x} } F(\mathbf{x} _0,{\bf 0} _N, \boldsymbol{ \theta })$, which is a $N \times N $ matrix that will be referred to as the {\bf connectivity matrix}, has a spectral radius smaller than one.
Since it is not possible to  compute the eigenvalues of $D_{\mathbf{x} } F(\mathbf{x} _0,{\bf 0} _N, \boldsymbol{ \theta }) $ for an arbitrary number of neurons $N$, we are hence obliged to proceed by finding estimations for the Cauchy bound~\cite{Rahman:Schmeisser} of its characteristic polynomial or by bounding the spectral radius $ \rho(D_{\mathbf{x} } F(\mathbf{x} _0,{\bf 0} _N, \boldsymbol{ \theta })) $ using either a matrix norm or the Gershgorin discs~\cite{horn:matrix:analysis}. An in-depth study of all these options showed that it is the use of the maximum row sum matrix norm $||| \cdot |||_{\infty}$ that yields the best stability bounds via the following statement:
\begin{equation}
\label{stability discrete main}
\rho(D_{\mathbf{x} } F(\mathbf{x} _0,{\bf 0} _N, \boldsymbol{ \theta }))\leq |||D_{\mathbf{x} } F(\mathbf{x} _0,{\bf 0} _N, \boldsymbol{ \theta })|||_{\infty}< 1 \quad \mbox{if and only if} \quad
|\partial _{{x}}f(x _0, 0, \boldsymbol{\theta}  )|<1.
\end{equation}
Notice that this remarkable result puts together the stability conditions for the continuous and discrete time systems.

As an example of application of these results, consider the Mackey-Glass nonlinear kernel~\cite{mackey-glass:paper}
\begin{equation}
\label{Mackey-Glass nonlinear kernel}
f(x,I, \boldsymbol{\theta})= \frac{\eta \left(x+\gamma I \right)}{1+ \left(x+\gamma I\right)^p},
\end{equation} 
where the parameter $\boldsymbol{\theta}:=(\gamma, \eta, p)$ is a three tuple of real values; $\gamma$  is usually referred to as the {\bf input gain} and $\eta $ the {\bf feedback gain}. When this prescription is used in~\eqref{time delay equation 2} in the autonomous regime, that is, $I (t)=0$, the associated dynamical system exhibits two families of equilibria $x _0$  parametrized by $\eta$, namely, $x _0=0 $ and the roots of $x_0 ^p= \eta-1  $. For example, in the case $p=2 $, two distinct cases arise: when $\eta< 1 $ there is a unique  equilibrium at the origin which is stable as long as $\eta\in (-1,1) $. When $\eta>1 $ two other equilibria appear at $x _0=\pm \left(\eta-1 \right)^{1/2}$ which are stable whenever $\eta<3 $. These statements are proved in Corollary~\ref{Stability MG continuous} of the Supplementary Material. Analogous statements for the Ikeda kernel~\cite{Ikeda1979}  
$
f(x,I, \boldsymbol{\theta})= \eta \sin ^2 \left(x+ \gamma I+ \phi\right)$,
$\boldsymbol{\theta}:=(\eta, \gamma, \phi)$
can be found in Corollary~\ref{suppl-Stability Ikeda continuous} of the Supplementary Material. A particularly convenient sufficient condition is $| \eta | \leq 1 $ that simultaneously ensures stability and unimodality (existence of a single stable equilibrium).

\paragraph{Empirical evidence.}

\begin{figure}[!ht]
\centering
\begin{tikzpicture}
    \node[anchor=south west,inner sep=0] at (0,0) {\hspace{-1cm}\includegraphics[scale=.37,trim = 3cm 1cm 1.5cm 1.5cm ,clip]{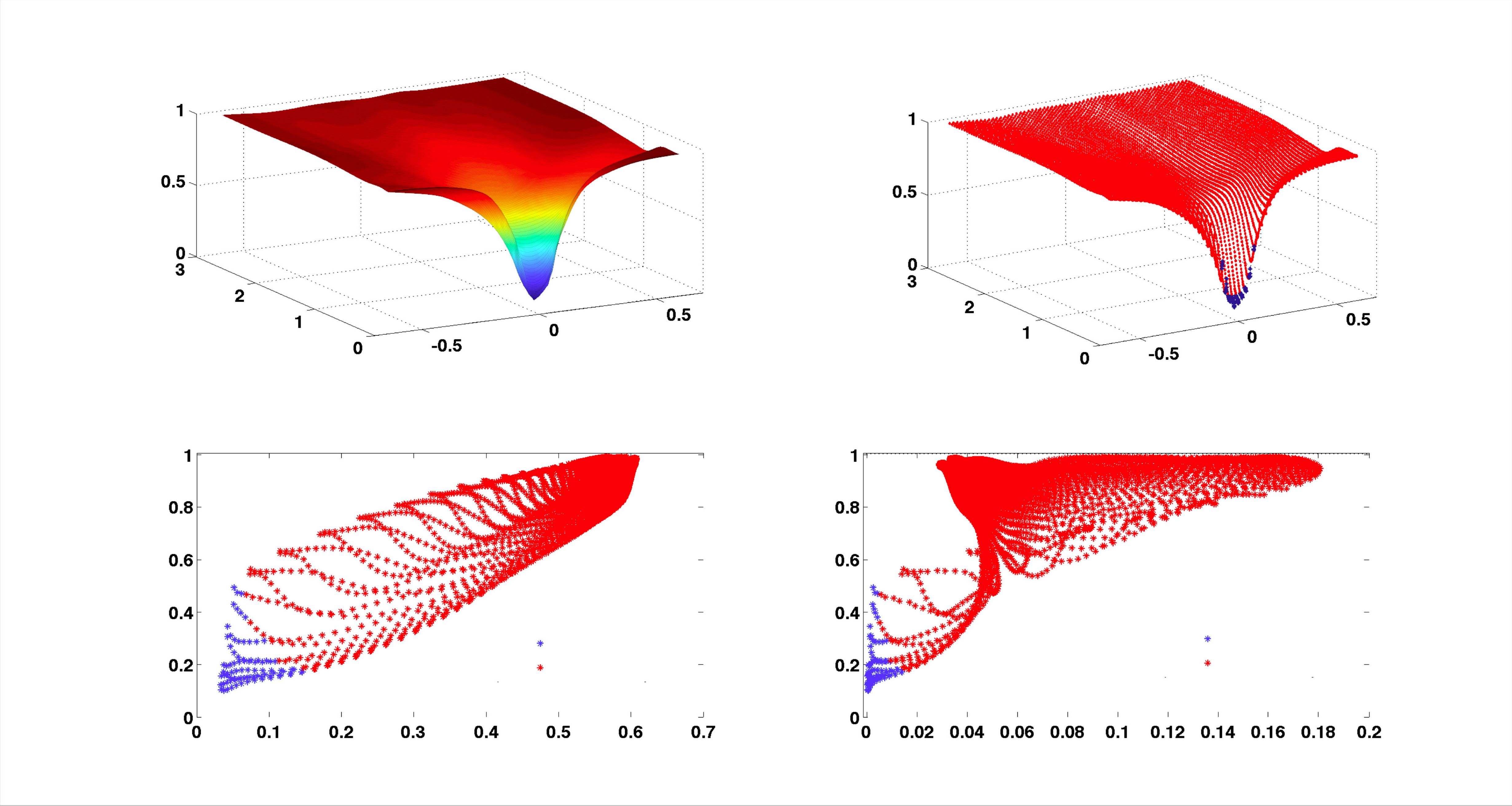}};
          \node at (5.51,1.63) {\tiny {\bf Unimodality}};          
          \node at (5.45,1.3) {\tiny {\bf Bimodality}};           
          \node at (13.75,1.66) {\tiny {\bf Unimodality}};    
          \node at (13.65,1.38)  {\tiny {\bf Bimodality}};
           \node at (12,0.2) {\scriptsize {\bf Variance of the reservoir output}}; 
          \node at (12,4.6) {\scriptsize {\bf Influence of the reservoir output variance on the }}; 
          \node at (12,4.3) {\scriptsize {\bf reservoir performance in the $3$-lag quadratic memory task}};  
            \node[label=below:\rotatebox{90}{\scriptsize {\bf NMSE}}] at (7.8,2.7) {};
            \node[label=below:\rotatebox{90}{\scriptsize {\bf NMSE}}] at (-0.3,2.7) {}; 
             \node at (3.7,0.2) {\scriptsize {\bf Mean of the reservoir output}}; 
          \node at (3.2,4.6) {\scriptsize {\bf Influence of the input mask mean on the reservoir }}; 
          \node at (3.2,4.3) {\scriptsize {\bf performance in the $3$-lag quadratic memory task}};
           \node[label=below:\rotatebox{90}{\scriptsize {\bf NMSE}}] at (8.6, 7.8) {};
           \node[label=below:\rotatebox{90}{\scriptsize {\bf NMSE}}] at (-0.39,7.8) {};
           \node at (12,9.3) {\scriptsize {\bf Influence of the input mask on the reservoir }}; 
          \node at (12,9) {\scriptsize  {\bf performance in the $3$-lag quadratic memory task}};
          \node at (3.5,9.3) {\scriptsize {\bf  Influence of the input mask on the reservoir }}; 
          \node at (3.5,9) {\scriptsize  {\bf performance in the $3$-lag quadratic memory task}};
          \node[label=below:\rotatebox{8}{\tiny {\bf Mean of the input mask}}] at (5.2,5.7) {};
          \node[label=below:\rotatebox{-24}{\tiny {\bf Variance of the input mask}}] at (0.7,6.7) {};
          \node[label=below:\rotatebox{10}{\tiny {\bf Mean of the input mask}}] at (13.6,5.7) {};
          \node[label=below:\rotatebox{-23}{\tiny {\bf Variance of the input mask}}] at (9.3,6.7) {};
\end{tikzpicture}
\caption{Behavior of the reservoir performance in a quadratic memory task as a function of the mean and the variance of the input mask. The modification of any of these two parameters influences how the reservoir dynamics separates from the stable equilibrium. The top panels show how the performance degrades very quickly as soon as the  mean and the variance of the input mask (and hence of the input signal) separate from zero. The bottom panels depict the reservoir performance as a function of the various output means and variances obtained when changing the input means and variances. In the top right and bottom panels we have indicated with  red markers the cases in which the reservoir visits the stability basin of a contiguous stable equilibrium hence showing how unimodality is associated to optimal performance.  }
\label{unimodality_via_simulatons}
\end{figure}
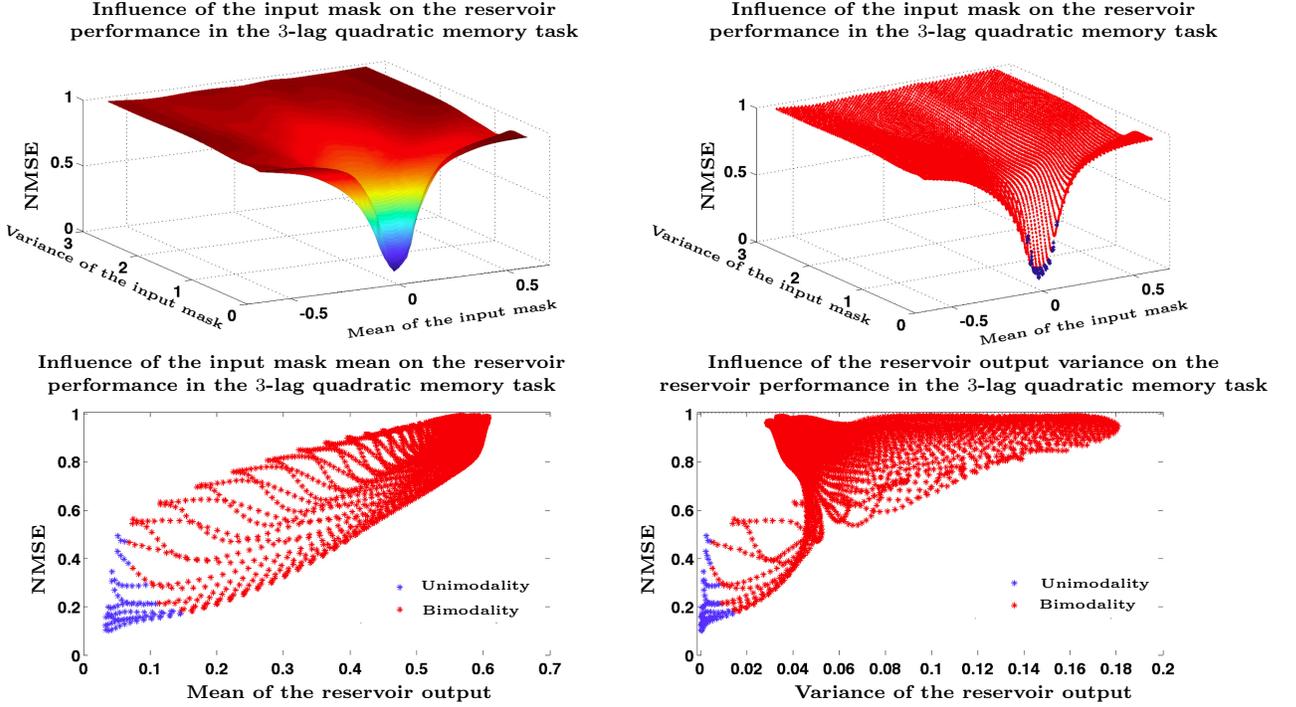

In order to confirm these theoretical and experimental arguments, we have carried out several numerical simulations in which we studied the RC performance in terms of the dynamical regime of the reservoir at the time of carrying out various nonlinear memory tasks. More specifically,  we construct a reservoir using the Ikeda nonlinear kernel  with $N = 20$, $d=0.2581$, $\eta=1.2443 $, $\gamma= 1.4762 $, and $\phi=0.1161 $. The equilibria of the associated autonomous system are given by the points $x _0$ where the curves $y=x$ and $y = \eta \sin ^2 \left(x+ \phi\right)$ intersect. With this parameter values, intersections take place at $x _0=0.0244 $, $x _0=0.9075$, and $x _0=1.063 $, which makes multi modality possible. As it can be shown with the results in the Supplementary Material section (see Corollary~\ref{suppl-Stability Ikeda continuous}), the first and the third equilibria are stable. 
In order to verify that the optimal performance is obtained  when the RC operates in the neighborhood of a stable equilibrium, we study the normalized mean square  error (NMSE) exhibited by a TDR initialized at $x _0=0.0244 $ when we present to it a quadratic memory task. More specifically, we inject in a TDR  an independent and identically normally distributed signal $z (t)$ with mean zero and variance $10^{-4} $ and we then train a linear readout ${\bf W}_{{\rm out}}$ (obtained with a ridge penalization of $\lambda=10^{-15}$) in order to recover the quadratic function $z (t-1)^2+  z (t-2)^2+z (t-3)^2$ out of the reservoir output. The top left panel in Figure~\ref{unimodality_via_simulatons} shows how the NMSE behaves as a function of the mean and the variance of the input mask $ {\bf c} $. It is clear that by modifying any of these two parameters we control how far the reservoir dynamics separates from the stable equilibrium, which we quantitatively evaluate in the two bottom panels by representing the RC performance in terms of the mean and the variance of the resulting reservoir output. Both panels depict how the injection of a signal slightly shifted in mean or with a sufficiently high variance results in reservoir outputs that separate from the stable equilibrium and in a severely degraded performance. An important factor in this deterioration seems to be the multi modality, that is, if the shifting in mean or the input signal variance are large enough then the reservoir output visits the stability basin of the other stable point placed at $x _0=1.063 $; in the top right and bottom panels we have marked with red color the values for which bimodality has occurred so that the negative effect of this phenomenon is noticeable. In the Supplementary Material section we illustrate how the behavior that we just described is robust with respect to the choice of nonlinear kernel and is similar when the experiment is carried out using the Mackey-Glass function.

\subsection{The approximating model and the nonlinear memory capacity of the reservoir computer}

The findings just presented have major consequences in the theoretical tools available for the evaluation of the RC performance. Indeed, since we now know that optimal operation is attained when the TDR functions in a unimodal fashion around an asymptotically stable steady state, we can approximate it by its partial linearization  with respect to the delayed self feedback term at that point and keeping the nonlinearity for the input injection. For statistically independent input signals of the type used to compute nonlinear memory capacities of the type introduced in~\cite{dambre2012}, this approximation allows us to visualize the TDR as a $N$-dimensional ($N$ is the number of neurons) vector autoregressive stochastic process of order one~\cite{luetkepohl:book} (we denote it as VAR(1)) for which the value of the associated nonlinear memory capacities can be explicitly computed. As we  elaborate later on in the discussion, the quality of this approximation at the time of evaluating the memory capacities of the original system is excellent and the resulting function can be hence used for RC optimization purposes regarding the nonlinear kernel parameter values $\boldsymbol{\theta} $ and the input mask ${\bf c} $.

Consider  a stable equilibrium $x _0 \in \mathbb{R}  $ of the autonomous system associated to~\eqref{time delay equation 2} or, equivalently, a stable fixed point of~\eqref{vector discretized reservoir main} of the form $\mathbf{x}_0:= (x _0, \ldots, x _0)^\top \in \mathbb{R} ^N$.  If we approximate~\eqref{vector discretized reservoir main} by its partial linearization at $\mathbf{x}_0 $ with respect to the delayed self feedback and by the $R$-order Taylor series expansion of the functional that describes the signal injection, we obtain an expression of the form:
\begin{equation}
\label{linearization with the random component main}
\mathbf{x} (t) = F(\mathbf{x} _0, {\bf 0}_N , \boldsymbol{\theta} ) +A(\mathbf{x} _0,  \boldsymbol{\theta} ) (\mathbf{x} (t-1) - \mathbf{x} _0)  + \boldsymbol{\varepsilon}  (t),
\end{equation}
where $A(\mathbf{x} _0,  \boldsymbol{\theta} ):=D_{\mathbf{x}} F(\mathbf{x} _0, {\bf 0}_N , \boldsymbol{\theta} )$ is the linear connectivity matrix and $\boldsymbol{\varepsilon} (t)$ is given by:
\begin{equation}
\label{eps_t_main}
\boldsymbol{ \varepsilon } (t) = (1 - { e} ^{- \xi})
\left(
  q _R \left( z(t), c _1\right)     ,
  q _R \left(z(t), c _1 , c _2 \right)    ,
  \ldots ,
  q _R \left(z(t), c _1 , \dots, c _N \right)
\right)^\top,
\end{equation}
with  
\begin{equation}
\label{q_R_polynomial_main}
q _R \left( z(t), c _1 , \dots, c _r \right) := \sum^{R}_{i = 1} \dfrac{z(t) ^i }{i !} (\partial _{{I} } ^{(i)} f   )({x} _0 , 0, \boldsymbol{\theta}  )  {\sum^{r}_{j = 1} { e} ^{-(r-j)\xi} c _j ^i} ,
\end{equation}
and $(\partial _{{I} } ^{(i)} f   )({x} _0 , 0, \boldsymbol{\theta}  )$ is the $i$th order partial derivative of the nonlinear kernel $f$ with respect to the second argument $I(t)$, evaluated at the point $({x} _0 , 0, \boldsymbol{\theta}  )$. 

If we now use as input signal  ${z} (t) $ independent and identically distributed random variables with mean $0$ and variance $\sigma _z ^2$ (we denote it by $\left\{ {z} (t)\right\}_{t \in \mathbb{Z} }\sim {\rm IID} (0, \sigma _z ^2 )$) then the recursion~\eqref{linearization with the random component main} makes the reservoir layer dynamics $\{ \mathbf{x} (t)\}_{t \in  \mathbb{Z} } $ into a discrete time random process that, as we show in what follows, is the solution of a $N$-dimensional vector autoregressive model of order 1 (VAR(1)). Indeed, it is easy to see that the assumption $\left\{ {z} (t)\right\}_{t \in  \mathbb{Z} }\sim {\rm IID} (0, \sigma _z ^2 )$ implies that $ \left\{ \mathbf{I} (t) \right\}_{t\in  \mathbb{Z}} \sim {\rm IID} ({\bf 0} _N , \Sigma _I  )$, with $\Sigma _I :=\sigma _z ^2 {\bf c}  ^\top {\bf c}$, and that $\{\boldsymbol{ \varepsilon } (t)\}_{t \in  \mathbb{Z}} $ is a family of $N$-dimensional independent and identically distributed random variables  with mean $\boldsymbol{\mu}_{{ \varepsilon }} $ and covariance matrix $\Sigma _{{ \varepsilon }} $ given by the following expressions:
\begin{equation}
\label{mu_eps_res_main}
\boldsymbol{\mu}_{{ \varepsilon }} = {\rm E} \left[ \boldsymbol{ \varepsilon } (t) \right] =  (1 - { e} ^{- \xi}) 
\left(
  q _R \left( \mu_{{z} }, c _1 \right)     ,
  q _R \left(\mu_{{z} }, c _1 , c _2 \right)    ,
  \ldots ,
  q _R \left(\mu_{{z} }, c _1 , \dots, c _N \right)
\right)^\top,
\end{equation}
where the polynomial $q _R$ is the same as in \eqref{q_R_polynomial_main} and where we use the convention that  the powers $ \mu_{{z} } ^i := {\rm E}\left[z(t)^i\right]$, for any $i \in \left\{ 1, \dots, R\right\} $ and with ${\rm E }[\cdot ]$ denoting the mathematical expectation. Additionally, $\Sigma _{{ \varepsilon }}:={\rm E}\left[ (\boldsymbol{ \varepsilon } (t)  - \boldsymbol{\mu}_{{ \varepsilon }}) (\boldsymbol{ \varepsilon } (t)  - \boldsymbol{ \mu}_{{ \varepsilon }}) ^\top \right]$  has entries determined by the relation:
\begin{equation*}
\label{Sigma_eps_main}
(\Sigma _{{ \varepsilon }}) _{ij} = (1 - { e} ^{- \xi }) ^2 ( ( q_R( \cdot, c _1 , \dots, c _i )\cdot q_R( \cdot, c _1 , \dots, c _j))(\mu_{{z} }) - q_R( \mu_{{z} }, c _1 , \dots, c _i )q_R( \mu_{{z} }, c _1 , \dots, c _j )),
\end{equation*}
where the first summand stands for the multiplication of the polynomials $q_R( \cdot, c _1 , \dots, c _i ) $ and $q_R( \cdot, c _1 , \dots, c _j )$ and the subsequent evaluation of the resulting polynomial at $\mu _z $, and the second one is made out of the multiplication of the evaluation of the two polynomials.

Using these observations, we can consider~\eqref{linearization with the random component main} as the prescription of a VAR(1) model driven by the independent noise $\{\boldsymbol{ \varepsilon } (t)\}_{t \in  \mathbb{Z}} $. If the nonlinear kernel $f$ satisfies the generic condition that the polynomial in $u$ given by $\det \left(\mathbb{I}_N-A(\mathbf{x} _0,  \boldsymbol{\theta} )u\right) $, does not have roots in and on the complex unit circle, then~\eqref{linearization with the random component main} has a second order stationary solution~\cite[Proposition 2.1]{luetkepohl:book} $\{ \mathbf{x} (t)\}_{t \in  \mathbb{Z} } $ with time-independent mean given by
\begin{equation}
\label{mu_x_main}
\boldsymbol{\mu} _{{x} } = {\rm E} \left[ \mathbf{x} (t)\right] = (I _N -  A(\mathbf{x} _0,  \boldsymbol{\theta} ) )^{-1} (F(\mathbf{x} _0 , {\bf 0} _N , \boldsymbol{\theta} ) - A(\mathbf{x} _0,  \boldsymbol{\theta} )  \mathbf{x} _0  + \boldsymbol{\mu} _{{\varepsilon}})
\end{equation}
and an also time independent autocovariance function $\Gamma (k):= {\rm E}\left[\left(\mathbf{x}(t)-\boldsymbol{\mu} _{{x} } \right) \left(\mathbf{x}(t-k)-\boldsymbol{\mu} _{{x} } \right)^\top\right]$, $k \in  \mathbb{Z} $, recursively determined the Yule-Walker equations (see~\cite{luetkepohl:book} for a detailed presentation). Indeed, $\Gamma(0)  $ is given by the vectorized equality:
\begin{equation}
\label{acvf0}
{\rm vec}(\Gamma(0) )= \left(\mathbb{I}_{N^2}-A(\mathbf{x} _0,  \boldsymbol{\theta} )\otimes A(\mathbf{x} _0,  \boldsymbol{\theta} )\right)^{-1} {\rm vec}(\Sigma _{{ \varepsilon }}),
\end{equation} 
which determines the higher order  autocovariances via the relation $\Gamma(k)=A(\mathbf{x} _0,  \boldsymbol{\theta} ) \Gamma(k-1) $ and the identity $\Gamma(-k)=\Gamma(k)^\top$. As we explain in the following paragraphs, the moments~\eqref{mu_eps_res_main},~\eqref{mu_x_main}, and~\eqref{acvf0} are all that is needed in order to characterize the memory capacities of the RC. 

A {\bf $h$-lag memory  task} is determined by a (in general nonlinear) function $H: \mathbb{R} ^{h+1}\rightarrow \mathbb{R}  $ that is used to generate a one-dimensional signal $y (t):=H(z (t), z (t-1), \ldots,z(t-h))$ out of the reservoir input. Given a TDR computer, the optimal linear readout ${\mathbf W}_{{\rm out}}$ adapted to the memory task $H$ is given by the solution of a ridge linear regression problem with regularization parameter $\lambda \in \mathbb{R} $ (usually tuned during the training phase via cross-validation) in which the covariates are the neuron values corresponding to the reservoir output and the explained variables are the values $\{y (t)\} $ of the memory task function. The  $H$-memory capacity $C_H( \boldsymbol{\theta}, {\bf c}, \lambda)$ of the TDR computer under consideration characterized by a nonlinear kernel $f$ with parameters $\boldsymbol{\theta}$, an input mask ${\bf c}  $, and a regularizing ridge parameter $\lambda$ is defined as one minus the normalized mean square error committed at the time of accomplishing the memory task $H$. When the reservoir is approximated by a VAR(1) process, then the corresponding {\bf $H$-memory capacity} is given by
\begin{equation}
\label{capacity formula}
C_H( \boldsymbol{\theta}, {\bf c}, \lambda)=\frac{{\rm Cov} (y (t) , \mathbf{x} (t))^\top(\Gamma(0) + \lambda \mathbb{I} _N ) ^{-1}(\Gamma (0) + 2\lambda \mathbb{I} _N) (\Gamma(0) + \lambda \mathbb{I} _N ) ^{-1}{\rm Cov} (y (t) , \mathbf{x} (t))}{{\rm var} \left(y (t)\right)}
\end{equation}
The developments leading to this expression are contained in the Supplementary Material section.
It is easy to show that:
\begin{equation*}
0 \leq C_H( \boldsymbol{\theta}, {\bf c}, \lambda)\leq 1.
\end{equation*}
Notice that in order to evaluate~\eqref{capacity formula} for a specific memory task, only ${\rm Cov} (y (t) , \mathbf{x} (t)) $ and ${\rm var} \left(y (t)\right) $ need to be computed since the autocovariance $\Gamma (0)  $ is fully determined by~\eqref{acvf0} once the reservoir and the equilibrium $\mathbf{x}_0 $ around which we operate have been chosen. As an example, we provide the expressions corresponding to the two most basic information processing routines, namely the linear and the quadratic memory tasks. Details on how to obtain the following equalities are contained in the Supplementary Material section.

\medskip

\noindent {\bf The $h$-lag linear memory task.} Linear memory tasks are those associated to linear task functions $H: \mathbb{R}^{h+1}\rightarrow \mathbb{R}$, that is, if we denote  ${\bf z}^h(t):= \left(z (t), z (t-1), \ldots,z(t-h)\right)^\top $ and $\mathbf{L}  \in \mathbb{R}^{h+1}$, we set $H({\bf z}^h(t)) :=\mathbf{L} ^\top {\bf z}^h(t) $. Various computations included in the Supplementary Material section using the so called  MA($\infty $) representation of the VAR(1) process show that 
${\rm var} ({y} (t) ) = \sigma _{{z} } ^2 ||\mathbf{L}  || ^2 $,
and
${\rm Cov} ({y} (t), x _i (t))= (1 - { e}^{- \xi }) \sum^{h+1}_{j = 1}  \sum^{N}_{s = 1}L_j (A(\mathbf{x} _0,  \boldsymbol{\theta} )^{j-1} )_{is} p_{R} (\mu_{z}, c _1 , \dots, c _s)$,  $i \in \{1, \ldots, N\} $,
where the polynomial $p_{R} $ on the variable $x$ is defined by  
$p_{R} (x, c _1 , \dots, c _s):=x \cdot q_{R} (x, c _1 , \dots, c _s)$ and its evaluation in the previous formula follows the same convention as in~\eqref{mu_eps_res_main}. 

\medskip

\noindent {\bf The $h$-lag quadratic memory task.} In this case we use a quadratic task function of the form 
\begin{equation}
\label{quadratic memory}
H({\bf z}^h(t)) := \mathbf{z} ^h (t) ^{\top} Q \mathbf{z} ^h (t) = \sum^{h+1}_{i = 1} \sum^{h+1}_{j = 1}Q _{ij} z({t-i+1}) z({t-j+1}),
\end{equation}
for some symmetric $h+1 $-dimensional matrix $Q$. If we define $y (t):=H({\bf z}^h(t)) $, we have that
$
{\rm var} ({y} (t) ) = (\mu _z ^4- \sigma _z ^4)\sum_{i=1}^{h+1} Q_{ii}^2+4\sigma_z ^4\sum_{i=1}^{h+1}\sum_{j>i}^{h+1}Q_{ij}^2$, 
and 
$${\rm Cov} ({y} (t), x _i(t))
 =  (1 - { e}^{- \xi }) \sum^{h+1}_{j = 1}  \sum^{N}_{r = 1} Q_{jj} (A ^{j-1} )_{ir}  (s_{R} (\mu_{{z} }, c _1 , \dots, c _r) - \sigma _{z}^2 q_R(\mu_z, c _1 , \dots, c _r)),$$
where the polynomial $s_R$ on the variable $x$ is defined as
$s_{R} (x, c _1 , \dots, c _r):=x ^2  \cdot q_{R} (x, c _1 , \dots, c _r)$. 

\section{Discussion}
\label{Discussion}

The possibility to approximate the TDR using a model of the type~\eqref{linearization with the random component main} opens the door to the theoretical treatment of many RC design related questions that so far were addressed using a trial and error approach. In particular,
the availability of a closed form formula of the type~\eqref{capacity formula} for the memory capacity of the RC is extremely convenient to determine the optimal reservoir architecture to carry out a given task. Nevertheless, it is obviously very important to assess the quality of the VAR(1) approximation underlying it and of the consequences that result from it. Indeed, we recall that the expression~\eqref{capacity formula}  was obtained via the partial linearization of the reservoir at a stable equilibrium in which it is initialized and kept in stationary operation. Despite the good theoretical and experimental reasons to proceed in this fashion explained in Section~\ref{Optimal performance: stability and unimodality}, we have confirmed their pertinence by explicitly comparing the reservoir memory capacity surfaces obtained empirically with those coming from the analytical expression~\eqref{capacity formula}. We have carried this comparison out for various tasks and have constructed the memory capacity surfaces as a function of different design parameters. 

\begin{figure}[!ht]
\hspace*{-2cm}\includegraphics[scale=.29]{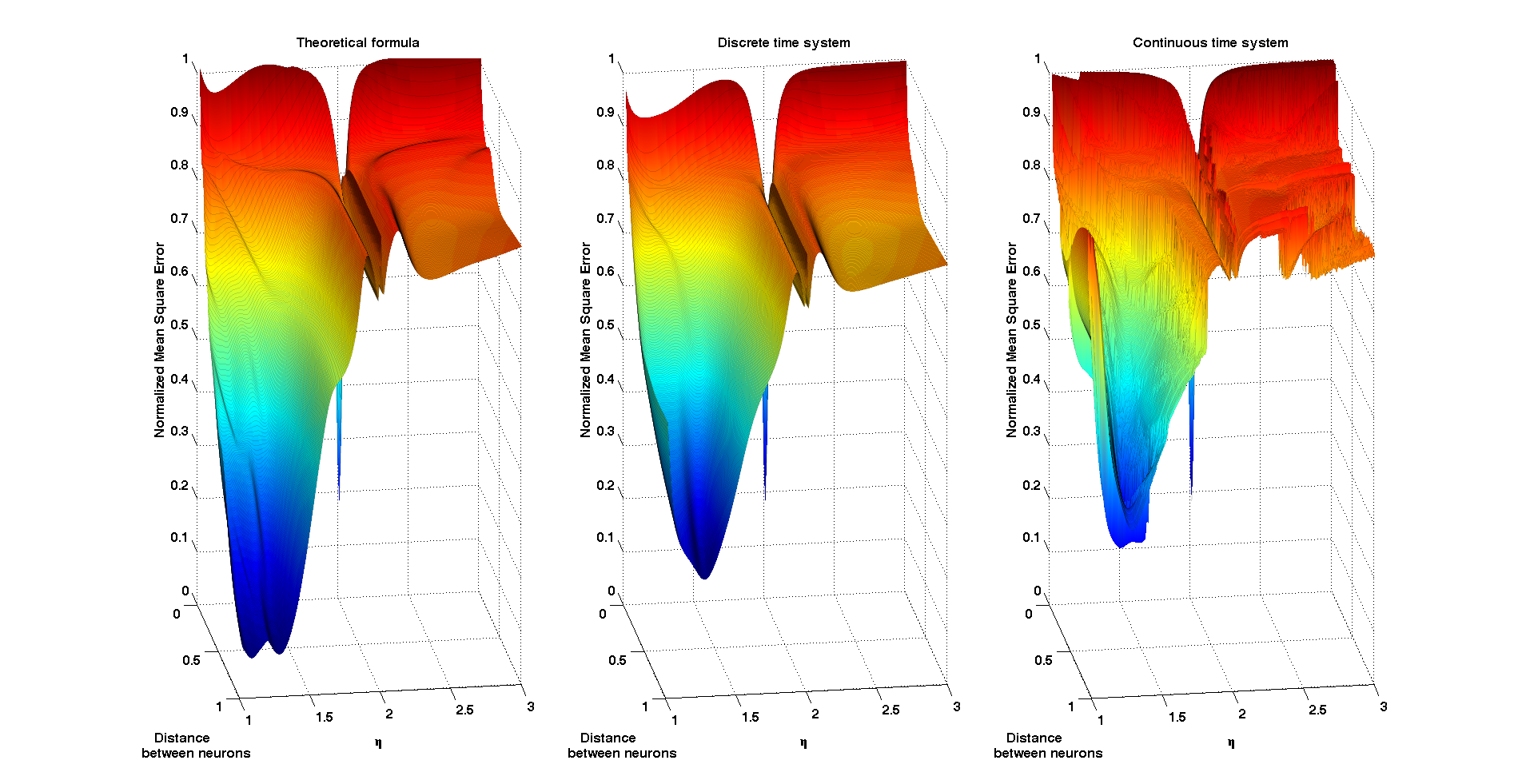}
\caption{Error surfaces exhibited by a Mackey-Glass kernel based reservoir computer in a 6-lag quadratic memory task, as a function of the distance between neurons and the parameter $\eta$.  The points in the surfaces of the middle and right panels are the result of Monte Carlo evaluations of the NMSE exhibited by the discrete and continuous time TDRs, respectively.  The left panel was constructed using the formula~\eqref{capacity formula} that is obtained as a result of modeling the reservoir with an approximating VAR(1) model.}
\label{memory as function eta theta} 
\end{figure}

We first consider a RC constructed using the Mackey-Glass nonlinear kernel~\eqref{Mackey-Glass nonlinear kernel} with $p=2$, $\gamma=0.796 $, and twenty neurons. We present to it the 6-lag quadratic memory task $H$ corresponding to choosing in~\eqref{quadratic memory} a seven dimensional diagonal matrix $Q$ with the diagonal entries given by the vector $(0,1,1,1,1,1,1) $. The first element, corresponding to the 0-lag memory (quadratic nowcasting), is set to zero in order to keep the difficulty of the task high enough. We then vary the value $d$ of the distance between neurons between 0 and 1 and the feedback  gain parameter $\eta $ between 1 and 3. As we discussed in Section~\ref{Optimal performance: stability and unimodality},  the TDR in autonomous regime exhibits for these parameter values two stable equilibria placed at $\pm (\eta-1)^{1/2}$; for this experiment we will always work with the positive equilibria by initializing the TDRs at those points. Figure~\ref{memory as function eta theta} represents the normalized mean square error (NMSE)  surfaces (which amounts to one minus the capacity) obtained using three different approaches. The left panel was obtained using the formula~\eqref{capacity formula} constructed a eight-order Taylor expansion of the nonlinear kernel on the signal input ($R=8 $ in~\eqref{eps_t_main}). The points in the surfaces of the middle and right panels are the result of Monte Carlo evaluations  (using 50,000 occurrences each) of the NMSE exhibited by the discrete and continuous time TDRs, respectively. The time-evolution of the time-delay differential equation (continuous time model) was simulated using a Runge-Kutta fourth-order method with a discretization step equal to $d/5 $.  A quick inspection of Figure~\ref{memory as function eta theta} reveals the ability of~\eqref{capacity formula} to accurately capture most of the details of the error surface and, most importantly, the location in parameter space where optimal performance is attained; it is very easy to visualize in this particular example how sensitive the magnitude of the error and the corresponding memory capacity are to the choice of parameters   and how small in size the region in parameter space associated with acceptable operation performance may be.

\begin{figure}[!ht]
\hspace*{-.1cm}\includegraphics[scale=.29]{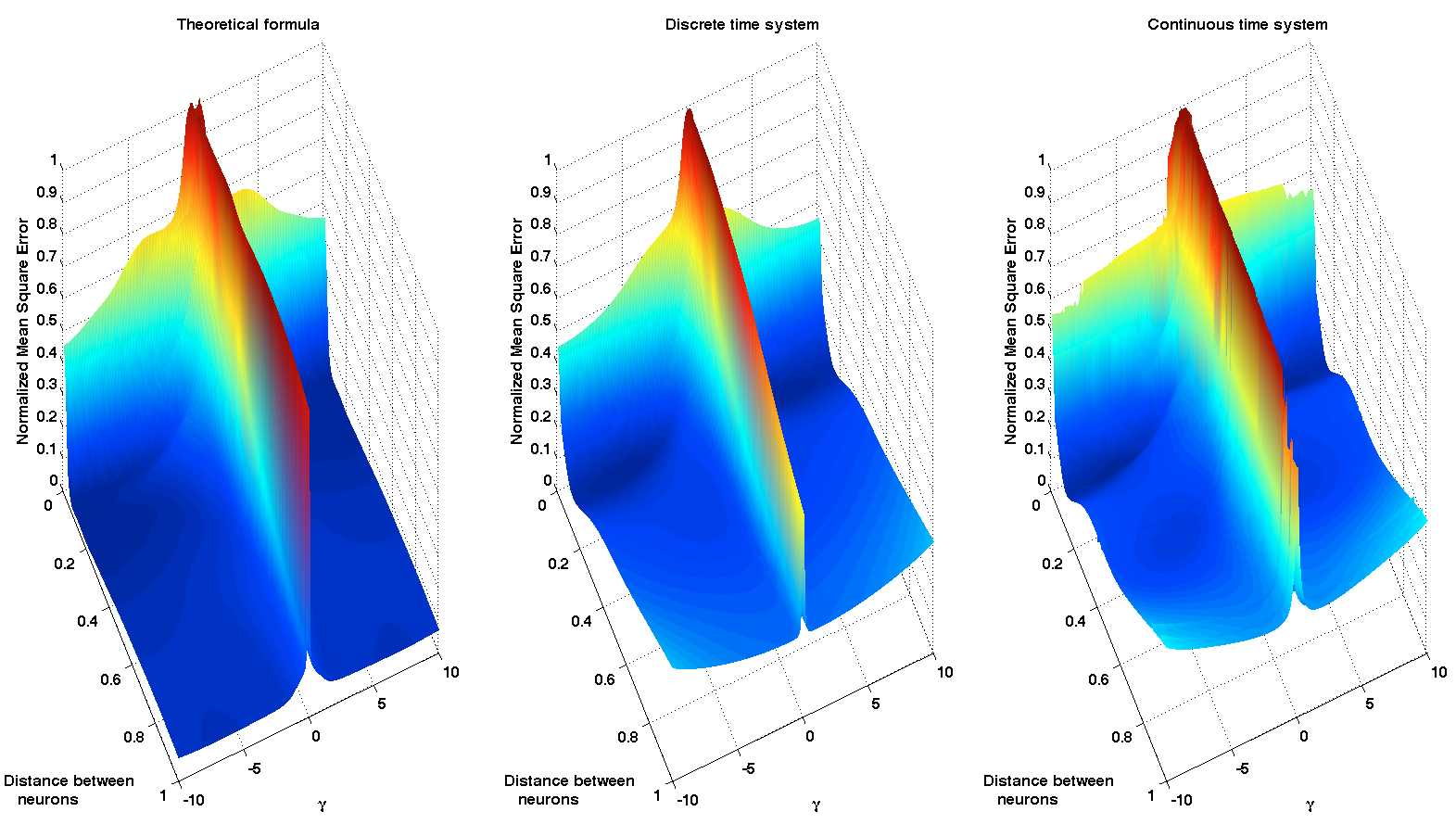}
\caption{Error surfaces exhibited by a Mackey-Glass kernel based reservoir computer in a 3-lag quadratic memory task, as a function of the distance between neurons and the parameter $\gamma$.  The points in the surfaces of the middle and right panels are the result of Monte Carlo evaluations of the NMSE exhibited by the discrete and continuous time TDRs, respectively.  The left panel was constructed using the formula~\eqref{capacity formula} that is obtained as a result of modeling the reservoir with an approximating VAR(1) model.}
\label{memory as function gamma theta}
\end{figure}

In order to show that these statements are robust with respect to the choice of task and varying parameters, we have carried out a similar experiment with a RC in which we fix the feedback gain $\eta _0=1.0781 $ and we vary the input gain $\gamma $ and the distance between neurons $d$. The quadratic memory task is reduced this time to 3-lags. We emphasize that in this setup the stable operation point is always the same and equal to $(\eta_0-1)^{1/2}$. Figure~\ref{memory as function gamma theta} shows how the performance of the memory capacity estimate~\eqref{capacity formula} at the time of capturing the optimal parameter region is in this situation comparable to the results obtained for the 6-lag quadratic memory task represented in Figure~\ref{memory as function eta theta}. We also point out that in this case there is a lower variability of the performance which, in our opinion, has to do with the fact that modifying the parameter $\gamma$ adjusts the input gain but leaves unchanged the operation point. Additionally, the moderate difficulty of the task makes possible attaining lower optimal error rates with the same number of neurons. In order to ensure the robustness of these results with respect to the choice of nonlinear kernel, we have included in the Supplementary Material section the results of a similar experiment carried out using the Ikeda prescription.

\begin{figure}[!ht]
\centering
\begin{tikzpicture}
    \node[anchor=south west,inner sep=0] at (0,0) {\hspace*{-2cm}\includegraphics[scale=.28]{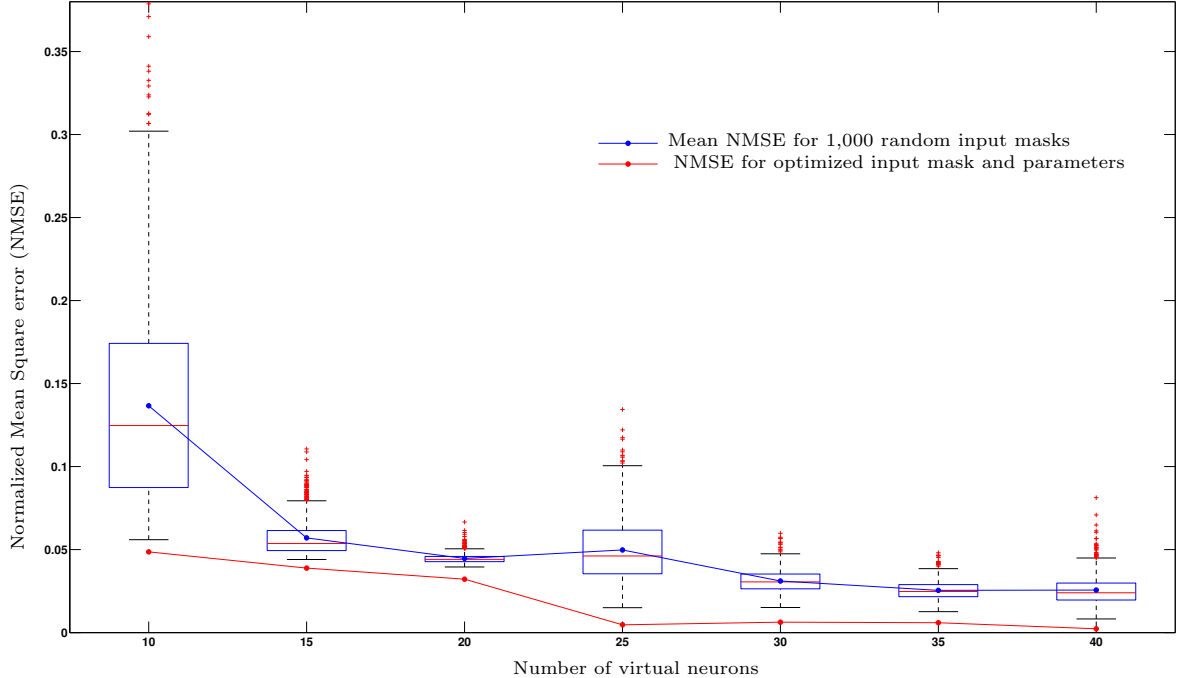}};                  
          \node at (11.1,7.8) {\scriptsize Mean NMSE for 1,000 random input masks}; 
          \node at (11.5,7.5) {\scriptsize NMSE for optimized input mask and parameters}; 
          \node at (8,0.8) {\scriptsize Number of virtual neurons}; 
           \node[label=below:\rotatebox{90}{\scriptsize Normalized Mean Square error (NMSE)}] at (-0.2,7.5) {};
\end{tikzpicture}
\caption{Influence of the mask optimization on the reservoir performance in the 3-lag quadratic memory task. The red line links the points that indicate the error committed by a RC with optimized parameters and mask. The box plots give information about the distribution of performances obtained with  1,000 input masks  randomly picked (only reservoir parameters have been optimized). As it is customary, on each box, the central mark is the median and the edges of the box are the 25th and 75th percentiles ($q _1$ and $q _3 $, respectively). The whiskers extend to the most extreme data points not considered outliers and outliers are plotted individually using red crosses. Points are drawn as outliers if they are larger than $q_3 + 1.5(q_3 Ð q_1) $ or smaller than $q_1-1.5(q_3 -q_1)$. The blue line links the points that indicate the mean NMSE committed when using the 1,000 different randomly picked masks.}
\label{influence of mask}
\end{figure}

Once the adequacy of the memory capacity evaluation formula~\eqref{capacity formula} has been established, we can use this result to investigate the influence of other architecture parameters in the reservoir performance. In Figure~\ref{influence of mask} we depict the results of an experiment where we study the influence of the choice of input mask $ {\bf c} $ in the performance of a Mackey-Glass kernel based reservoir in a 3-lag quadratic memory task. The figure shows, for each number of neurons, the performance obtained by a RC in which the reservoir parameters $\boldsymbol{\theta} $  and the input mask ${\bf c}  $ have been chosen so that the memory capacity $C_H( \boldsymbol{\theta}, {\bf c}, \lambda) $ in~\eqref{capacity formula} is maximized; we have subsequently kept the optimal parameters $\boldsymbol{\theta} $ and we have randomly constructed one thousand input masks ${\bf c}$  with entries belonging to the interval $[-3,3]$. The box plots in Figure~\ref{influence of mask} give an idea of the distribution of the degraded performances with respect to the optimal mask for different numbers of virtual neurons.

In conclusion,  the construction of  approximating models for the reservoir as well as the availability of performance evaluation formulas like~\eqref{capacity formula} based on it, constitute extremely valuable analytical tools whose existence should prove very beneficial in the fast and efficient extension and customization of RC type techniques to tasks far more sophisticated than the ones we considered in this paper. This specific point is the subject of ongoing research on which we will report in a forthcoming publication. 

\medskip
 
\noindent {\bf Acknowledgments:} We acknowledge partial financial support of the R\'egion de Franche-Comt\'e (Convention 2013C-5493), the European project PHOCUS (FP7 Grant
No. 240763), the Labex ACTION program (Contract No. ANR-11-LABX-01-01), and Deployment S.L. LG acknowledges financial support from the Faculty for the Future Program of the Schlumberger Foundation.
\newpage
\begin{center}
\Large{\textbf{Supplementary material for the paper ``Optimal nonlinear information processing capacity in delay-based reservoir computers''}}
\end{center}

\begin{appendices}
    \numberwithin{equation}{section}
    \numberwithin{figure}{section}
\section{Notation}
\label{Notation}

Column vectors are denoted by bold lower or upper case  symbol like $\mathbf{v}$ or $\mathbf{V}$. We write $\mathbf{v} ^\top $ to indicate the transpose of $\mathbf{v} $. Given a vector $\mathbf{v} \in \mathbb{R}  ^n $, we denote its entries by $v_i$, with $i \in \left\{ 1, \dots, n
\right\} $; we also write $\mathbf{v}=(v _i)_{i \in \left\{ 1, \dots, n\right\} }$. 
The symbols $\mathbf{i} _n$ and $ \mathbf{0} _n $ stand for the vectors of length $n$ consisting of zeros and ones, respectively. 

We denote by $\mathbb{M}_{n ,  m }$ the space of real $n\times m$ matrices with $m, n \in \mathbb{N} $. When $n=m$, we use the symbol $\mathbb{M}  _n $  to refer to the space of square matrices of order 
$n$. Given a matrix $A \in \mathbb{M}  _{n , m} $, we denote its components by $A _{ij} $ and we write $A=(A_{ij})$, with $i \in \left\{ 1, \dots, n\right\} $, $j \in \left\{ 1, \dots m\right\} $.   
We use $\mathbb{S}_n   $ to  indicate the subspace $\mathbb{S}  _n \subset \mathbb{M}  _n $ of symmetric matrices, that is, $\mathbb{S}  _n = \left\{ A \in \mathbb{M}  _n \mid A ^\top = A\right\}$.
Given a matrix $A \in \mathbb{M} _{n, m}$, the maximum row sum matrix norm  is defined as $|||A|||_{\infty} = \max_{1 \le i \le n} \sum^{m}_{j = 1} |A_{ij}|$. The symbol $\otimes$ stands for the Kronecker matrix product.

We use  ${C}^r([a,b], \mathbb{R})$, $r\ge 0$, to denote the Banach space of the $r$-times continuously differentiable real valued maps defined on the interval $[a,b]$
 with the topology of the uniform convergence; given a function $g \in {C}^r([a,b], \mathbb{R})$ when $r \ge 1$. We  designate the $l_{\infty}$-norm of an element $\phi \in {C}([a ,b], \mathbb{R})$ by $||\phi  || _\infty  = \sup_{\theta \in \left[ a , b\right] } | \phi(\theta ) |$.

The symbols ${\rm E}[\cdot]$, ${\rm var}(\cdot)$, and ${\rm Cov}(\cdot,\cdot)$ denote the mathematical expectation, the variance, and the covariance, respectively.
 \section{The reservoir map and the connectivity matrix}
\label{Preliminaries}

The reservoir map $F: \mathbb{R} ^N \times\mathbb{R}  ^N \times\mathbb{R}  ^K \longrightarrow \mathbb{R} ^N$ introduced in~\eqref{vector discretized reservoir main} is uniquely determined by the recursions~\eqref{recursion euler} obtained out of the Euler discretization of the time-delay differential equation (TDDE)
\begin{equation}
\label{time delay equation}
 \dot{x} (t)= -x(t) + f(x (t- \tau), I(t), \boldsymbol{\theta}),
\end{equation}
and organized in neuron layers parametrized by $t \in  \mathbb{Z}$. The reservoir map is obtained by using~\eqref{recursion euler} in order to write down the neuron values of the layer for time $t$ in terms of those for time $t-1 $ and the current input signal value. More specifically:
\begin{equation}
\label{discretized reservoir}
\left\{
\begin{array}{rcl}
x_1 (t) &=& e  ^{- \xi } x _N (t-1) + (1- e ^{- \xi }) f(x _1 (t-1), I _1 (t), \boldsymbol{ \theta }),\\
x_2 (t) &=& e  ^{-2 \xi } x _N (t-1) + (1- e ^{- \xi }) \left\{ e  ^{- \xi }f(x _1 (t-1), I _1 (t), \boldsymbol{ \theta }) + f(x _2 (t-1), I _2 (t), \boldsymbol{ \theta })\right\},\\
  &\vdots& \\
x_N (t) &=& e  ^{-N \xi } x _N (t-1) + (1- e ^{- \xi })  \sum^{N-1}_{j = 0} e  ^{- j\xi }f(x _{N-j} (t-1), I _{N-j} (t), \boldsymbol{ \theta }),
\end{array}
\right.
\end{equation}
which corresponds to a description of the form
\begin{equation}
\label{vector discretized reservoir}
\mathbf{x}(t) = F( \mathbf{x} (t-1), \mathbf{I} (t), \boldsymbol{ \theta }) , 
\end{equation}
that uniquely determines the reservoir map $F: \mathbb{R} ^N \times\mathbb{R}  ^N \times\mathbb{R}  ^K \longrightarrow \mathbb{R} ^N$.

Let  $x _0 \in \mathbb{R}   $  and $\mathbf{x} _0:= x _0 \mathbf{i} _N \in \mathbb{R}^N  $. Let  $D_{\mathbf{x} } F(\mathbf{x} _0,{\bf 0} _N, \boldsymbol{ \theta }) $ be the partial derivative of $F$ with respect to the first argument computed at the point $(\mathbf{x} _0,{\bf 0} _N, \boldsymbol{ \theta })$. We will refer to $A(\mathbf{x} _0, \boldsymbol{ \theta }):=D_{\mathbf{x} } F(\mathbf{x} _0,{\bf 0} _N, \boldsymbol{ \theta }) $ as the {\bf  connectivity matrix} of the reservoir at the point $\mathbf{x} _0 $. It is easy to show that $A(\mathbf{x} _0, \boldsymbol{ \theta })$ has the following explicit form
\begin{align}
\label{DxF}
&A(\mathbf{x} _0, \boldsymbol{ \theta }):=D_{\mathbf{x} } F(\mathbf{x} _0,{\bf 0} _N, \boldsymbol{ \theta }) =  
{\small\left(
\begin{array}{cccccc}
 \Phi&    0 &\hdots &0&  { e^{- \xi }}{}\\
 { e^{- \xi }}\Phi&   \Phi&  \hdots &0& { e^{-2 \xi }}{}\\
  e^{- 2\xi }\Phi&  e^{- \xi }\Phi&  \hdots &0&  { e^{- 3\xi }}{}\\
  \vdots&\vdots&\ddots&\vdots&\vdots\\
  e^{- (N-1)\xi }\Phi&   e^{- (N-2)\xi }\Phi&   \hdots &{ e^{- \xi }}{} \Phi &\Phi + { e^{- N\xi }}{}
\end{array}
\right)},
\end{align}
where $ \Phi :=(1- e ^{- \xi })\partial _{x}f(x _0, 0 , \boldsymbol{\theta}  )$ and $\partial _{x}f(x _0, 0 , \boldsymbol{\theta}  ) $ is the first derivative of the nonlinear kernel $f$ in \eqref{time delay equation} with respect to the first argument and computed at the point $(x _0, 0 , \boldsymbol{\theta}  )$. We will also use  the symbol $f'_{x_0}$ to denote $\partial _{x}f(x _0, 0 , \boldsymbol{\theta})$.

\section{The approximating vector autoregressive system and information processing capacity estimations}
\label{VAR representation}
The goal of this appendix is providing details on the construction of the approximating VAR(1) model for the TDR obtained after a partial linearization of the reservoir map on the dynamical variables and respecting the nonlinearity on the input signal.  Let  $x _0 \in \mathbb{R}   $  and $\mathbf{x} _0:= x _0 \mathbf{i} _N \in \mathbb{R}^N  $ be an equilibrium of~\eqref{time delay equation} and a fixed point of~\eqref{vector discretized reservoir}, respectively (see Proposition~\ref{fixed point vs equilibrium}). These solutions are chosen with respect to the autonomous regime, that is, we set  $I(t) = 0$ in the time-delay differential equation \eqref{time delay equation} and $\mathbf{I} (t) = {\bf 0}_N $ in the associated recursion \eqref{vector discretized reservoir}. In practice, we choose stable solutions; Appendix~\ref{Equilibria of time-delay differential equation and those of its discretization} contains conditions that ensure this dynamical feature.

\medskip

\noindent{\bf The VAR model setup and the construction of the approximating VAR system.} The vector autoregressive (VAR) model is of much use in multivariate time series analysis and is a natural extension of the univariate linear autoregressive (AR) model. See~\cite{luetkepohl:book} for an extensive introduction to multivariate time series methods and the details on the VAR processes. The $N$-dimensional VAR(1) process (VAR model of order 1)  is defined in mean-adjusted form as the solution to the recursions
\begin{equation}
\label{VAR model}
\mathbf{x} (t) -   \boldsymbol{\mu}_x =  A  (\mathbf{x} (t-1) - \boldsymbol{\mu}_x) + \boldsymbol{ \epsilon }(t), \enspace t=0, \pm 1, \pm 2, \ldots
\end{equation}
where $\mathbf{x} (t) = (x _1 (t), \ldots, x _N (t))^\top \in \mathbb{R} ^N $ is a random vector, $A\in \mathbb{M} _N $ is a fixed coefficient matrix, $\boldsymbol{\mu}_x \in \mathbb{R} ^N $, and $\boldsymbol{ \epsilon }(t) = (\epsilon _1 (t), \ldots, \epsilon _N (t))^\top\in \mathbb{R} ^N $ is such that $\left\{\boldsymbol{ \epsilon }(t)\right\} \sim {\rm WN}({\bf 0} _N, \Sigma _{\epsilon})$ is a $N$-dimensional white noise (stochastic process that presents no autocorrelation) or innovation process with mean ${\bf 0} _N$ and  covariance matrix $\Sigma _{ \epsilon } \in  \mathbb{S} _N $. We are particularly interested in {\bf stable} VAR models, that is, models of the type~\eqref{VAR model} where the autoregression matrix $A$ is chosen such that 
\begin{equation}
\label{stationarity condition}
\det \left(\mathbb{I} _N-Az\right)\neq 0 \quad \mbox{for all $z \in \mathbb{C} $ such that $|z| \leq 1$}.
\end{equation}
It can be proved (see Proposition 2.1 in~\cite{luetkepohl:book}) that stable models have a unique second order stationary solution $\{ \mathbf{x}(t)\} _{t \in   \mathbb{Z}}$ for which $\boldsymbol{\mu}_x={\rm E}[\mathbf{x}(t)] $ and the autocovariance function $$\Gamma (k):= {\rm E}\left[\left(\mathbf{x}(t)-\boldsymbol{\mu} _{{x} } \right) \left(\mathbf{x}(t-k)-\boldsymbol{\mu} _{{x} } \right)^\top\right], \enspace k \in   \mathbb{Z} $$ is time independent.

Let now $ \mathbf{x} _0=x _0 \mathbf{i}_N \in \mathbb{R} ^N $ be the stable fixed point of the reservoir map \eqref{vector discretized reservoir} in  autonomous regime, that is,   $F( {\bf x} _0  , {\bf 0} _N , \boldsymbol{ \theta }) = {\bf x} _0 $. We now approximate \eqref{vector discretized reservoir} by its partial linearization at $\mathbf{x} _0 $ with respect to the delayed self feedback and by the $R$th-order Taylor series expansion on the variable that determines the input signal injection. We obtain the following expression: 
\begin{equation}
\label{linearization of the reservoir}
\mathbf{x} (t) = F(\mathbf{x} _0, {\bf 0}_N, \boldsymbol{ \theta } ) + D_{\mathbf{x} } F(\mathbf{x} _0,{\bf 0} _N, \boldsymbol{ \theta }) (\mathbf{x} (t-1) - \mathbf{x} _0)  + \boldsymbol{\varepsilon}  (t),
\end{equation}
where  $D_{\mathbf{x} } F(\mathbf{x} _0,{\bf 0} _N, \boldsymbol{ \theta }) $ is the first derivative of $F$ with respect to its first argument, computed at the point $(\mathbf{x} _0,{\bf 0} _N, \boldsymbol{ \theta })$. We recall that $D_{\mathbf{x} } F(\mathbf{x} _0,{\bf 0} _N, \boldsymbol{ \theta })=A(\mathbf{x} _0, \boldsymbol{ \theta }) $ is the connectivity matrix introduced in~\eqref{DxF}. Additionally, $\boldsymbol{\varepsilon} (t)$ in \eqref{linearization of the reservoir} is obtained out of the Taylor series expansion of $F(\mathbf{x} (t), \mathbf{I} (t), \boldsymbol{\theta} )$ in \eqref{discretized reservoir} on $\mathbf{I} (t)$ up to some fixed order $R\in \mathbb{N} $ and is given by
\begin{equation}
\label{eps_t}
\boldsymbol{ \varepsilon } (t) = (1 - { e} ^{- \xi})(q _R \left( z(t), c _1\right), q _R \left(z(t), c _1 , c _2 \right), \dots, q _R \left(z(t), c _1 , \dots, c _N \right))^\top,
\end{equation}
with
\begin{equation}
\label{q_R_polynomial}
q _R \left( z(t), c _1 , \dots, c _r \right) := \sum^{R}_{i = 1} \dfrac{z(t) ^i }{i !} (\partial _{{I} } ^{(i)} f   )({x} _0 , 0, \boldsymbol{\theta}  )  {\sum^{r}_{j = 1} { e} ^{-(r-j)\xi} c _j ^i} ,
\end{equation}
where $c _i $, $i \in \left\{ 1, \ldots, N \right\} $ are the entries of the input mask ${\bf c} \in \mathbb{R} ^N $ and $(\partial _{{I} } ^{(i)} f   )({x} _0 , 0, \boldsymbol{\theta}  )$ is the $i$th order partial derivative of the nonlinear reservoir kernel $f$ in \eqref{time delay equation} with respect to the second argument $I(t)$ computed at the point $({x} _0 , 0, \boldsymbol{\theta}  )$.

If we now use as input signal  ${z} (t) $ independent and identically distributed random variables with mean $0$ and variance $\sigma _z ^2$, that is, $\left\{ {z} (t)\right\}_{t \in  \mathbb{Z} }\sim {\rm IID} (0, \sigma _z ^2 )$, then the recursion~\eqref{linearization of the reservoir} makes the reservoir layer dynamics $\{ \mathbf{x} (t)\}_{t \in  \mathbb{Z} } $ into a discrete time random process that, as we show in what follows, is the solution of a $N$-dimensional VAR(1) model. Indeed, it is easy to see that the assumption $\left\{ {z} (t)\right\}_{t \in  \mathbb{Z} }\sim {\rm IID} (0, \sigma _z ^2 )$ implies that $ \left\{ \mathbf{I} (t) \right\}_{t\in  \mathbb{Z}} \sim {\rm IID} ({\bf 0} _N , \Sigma _I  )$, with $\Sigma _I :=\sigma _z ^2 {\bf c}  ^\top {\bf c}$, and that $\{\boldsymbol{ \varepsilon } (t)\}_{t \in  \mathbb{Z}} $ is a family of $N$-dimensional independent and identically distributed random variables  with mean $\boldsymbol{\mu}_{{ \varepsilon }} $ and covariance matrix $\Sigma _{{ \varepsilon }} $ given by the following expression:
\begin{equation}
\label{mu_eps_res_main_suppl}
\boldsymbol{\mu}_{{ \varepsilon }} = {\rm E} \left[ \boldsymbol{ \varepsilon } (t) \right] =  (1 - { e} ^{- \xi}) 
\left(
  q _R \left( \mu_{{z} }, c _1 \right)     ,
  q _R \left(\mu_{{z} }, c _1 , c _2 \right)    ,
  \ldots ,
  q _R \left(\mu_{{z} }, c _1 , \dots, c _N \right)
\right)^\top,
\end{equation}
where the polynomial $q _R$ is the same as in \eqref{q_R_polynomial} and where we use the convention that  the powers $ \mu_{{z} } ^i := {\rm E}\left[z(t)^i\right]$, for any $i \in \left\{ 1, \dots, R\right\} $. For example, if  the variables $z (t)$ are normal, that is, $\left\{ {z} (t)\right\}_{t \in  \mathbb{Z} }\sim {\rm IN} (0, \sigma _z ^2 )$, then
\begin{equation}
\label{mu_I_res}
\mu_{{z} } ^i := {\rm E}\left[z(t)^i\right] = \left\{ \begin{array}{cl}
  \dfrac{2l!}{2^l l!} \sigma_z ^{2l}&  {\rm when} \enspace i = 2l,\quad l \in \mathbb{N},  \\
  & \\
  0&{\rm otherwise}.
\end{array}
\right.
\end{equation}
Additionally, $\Sigma _{{ \varepsilon }}:={\rm E}\left[ (\boldsymbol{ \varepsilon } (t)  - \boldsymbol{\mu}_{{ \varepsilon }}) (\boldsymbol{ \varepsilon } (t)  - \boldsymbol{ \mu}_{{ \varepsilon }}) ^\top \right]$  has entries determined by the relation:
\begin{equation*}
\label{Sigma_eps_main_suppl}
(\Sigma _{{ \varepsilon }}) _{ij} = (1 - { e} ^{- \xi }) ^2 ( ( q_R( \cdot, c _1 , \dots, c _i )\cdot q_R( \cdot, c _1 , \dots, c _j))(\mu_{{z} }) - q_R( \mu_{{z} }, c _1 , \dots, c _i )q_R( \mu_{{z} }, c _1 , \dots, c _j )),
\end{equation*}
where the first summand stands for the multiplication of the polynomials $q_R( \cdot, c _1 , \dots, c _i ) $ and $q_R( \cdot, c _1 , \dots, c _j )$ and the subsequent evaluation of the resulting polynomial at $\mu _z $, and the second one is made out of the multiplication of the evaluation of the two polynomials.

With these observations it is clear that we can consider~\eqref{linearization of the reservoir} as a VAR(1) model driven by the independent noise $\{\boldsymbol{ \varepsilon } (t)\}_{t \in  \mathbb{Z}} $. If the nonlinear kernel $f$ satisfies the generic condition that the polynomial in $z$ given by $\det \left(\mathbb{I}_N-A(\mathbf{x} _0,  \boldsymbol{\theta} )z\right) $, does not have roots in and on the complex unit circle, then~\eqref{linearization of the reservoir} has a unique second order stationary solution $\{ \mathbf{x} (t)\}_{t \in  \mathbb{Z} } $ with time-independent mean 
\begin{equation}
\label{mu_x}
\boldsymbol{\mu} _{{x} } = {\rm E} \left[ \mathbf{x} (t)\right] = (I _N -  A(\mathbf{x} _0,  \boldsymbol{\theta} ) )^{-1} (F(\mathbf{x} _0 , {\bf 0} _N , \boldsymbol{\theta} ) - A(\mathbf{x} _0,  \boldsymbol{\theta} )  \mathbf{x} _0  + \boldsymbol{\mu} _{{\varepsilon}}).
\end{equation}
that can be used to rewrite \eqref{linearization of the reservoir} in mean-adjusted form
\begin{equation}
\label{VAR model reservoir}
\mathbf{x} (t) - \boldsymbol{\mu} _{{x} } = A(\mathbf{x} _0,  \boldsymbol{\theta} )( \mathbf{x} (t-1) - \boldsymbol{\mu}  _{{x} }) + (\boldsymbol{\varepsilon} (t) - \boldsymbol{\mu}_{{\varepsilon}}).
\end{equation}
In the presence of stationarity we can recursively compute the time independent autocovariance function $\Gamma (k):= {\rm E}\left[\left(\mathbf{x}(t)-\boldsymbol{\mu} _{{x} } \right) \left(\mathbf{x}(t-k)-\boldsymbol{\mu} _{{x} } \right)^\top\right]$ at lag $k \in  \mathbb{Z} $ by  using the Yule-Walker equations~\cite{luetkepohl:book}. Indeed, $\Gamma(0)  $ is given by the vectorized equality:
\begin{equation}
\label{Gamma0} 
{\rm vec}(\Gamma(0) )= \left(\mathbb{I}_{N^2}-A(\mathbf{x} _0,  \boldsymbol{\theta} )\otimes A(\mathbf{x} _0,  \boldsymbol{\theta} )\right)^{-1} {\rm vec}(\Sigma _{{ \varepsilon }}),
\end{equation} 
which determines the higher order autocovariances  via the relation 
\begin{equation}
\label{Gammah}
\Gamma(k)=A(\mathbf{x} _0,  \boldsymbol{\theta} ) \Gamma(k-1),
\end{equation}
and the identity $\Gamma(-k)= \Gamma(k)^\top $.

\medskip

\noindent {\bf The nonlinear memory capacity estimations.}
We now concentrate on the computation of the quantitative measures of the reservoir performance introduced in the paper. In particular, we will provide details on the computation of the nonlinear memory capacity formula in~\eqref{capacity formula}.  Recall that a {\bf $h$-lag memory  task} is determined by a  function $H: \mathbb{R} ^{h+1}\rightarrow \mathbb{R}  $ (in general nonlinear) that is used to generate a one-dimensional signal $y (t):=H(z (t), z (t-1), \ldots,z(t-h))$ out of the reservoir input $\left\{z(t)\right\} _{t \in \mathbb{Z} }$. 

Consider now a TDR computer with $N$ neurons. The optimal linear readout ${\mathbf W}_{{\rm out}}$ adapted to the memory task $H$ is given by the solution of a ridge (or Tikhonov~\cite{tikhonov:regression}) linear regression problem with regularization parameter $\lambda \in \mathbb{R} $ (usually tuned during the training phase via cross-validation) in which the covariates are the neuron values corresponding to the reservoir output and the explained variables are the values $\{y (t)\} $ of the memory task function. More explicitly, ${\mathbf W}_{{\rm out}}$   is given by the solution of the following optimization problem
\begin{equation}
({\mathbf W}_{{\rm out}}, a_{\rm out}):=\mathop{\rm arg\, min}_{{\mathbf W}   \in \mathbb{R} ^N,  a\in\mathbb{R}  } \left({\rm E} \left[  ({\mathbf W} ^{ \top} \cdot {\bf x}(t)  +a -{ y}(t) )^2\right]  + \lambda \|{\mathbf W}   \|^2\right),\label{RC optimization problem}
\end{equation}
where the expectation is taken thinking of $y _t $ and $ \mathbf{x} (t)$ as random variables due to the stochastic nature of the input signal $\left\{  z(t) \right\} _{t\in  \mathbb{Z} }$ and hence that of the $\left\{ \mathbf{I} (t)\right\}_{t\in  \mathbb{Z} } $. 
In order to obtain the explicit solution of \eqref{RC optimization problem}, we first define $g({\mathbf W}, a):={\rm E} \left[  ({\mathbf W} ^{ \top} \cdot {\bf x}(t)  +a -{ y}(t) )^2\right]  + \lambda \|{\mathbf W}  \|^2 $ and set  
\begin{align*}
\dfrac{\partial g({\mathbf W},a)}{\partial w _i } = &2 \left[\sum^{N}_{j = 1} w _j {\rm Cov} (x_j (t), x _i (t)) - {\rm Cov} (y (t) , x _i (t)) + \lambda w_i \right] = 0, \enspace i\in \left\{ 1, \ldots, N \right\} ,\\
\dfrac{\partial g({\mathbf W},a)}{\partial a} = &2 \left[ a + {\mathbf W} ^\top {\rm E} \left[ \mathbf{x} (t)\right]  - {\rm E} \left[ y (t) \right] \right] = 0,
\end{align*}
or, equivalently,
\begin{align*}
({\rm Cov}(\mathbf{x} (t), \mathbf{x} (t)) + \lambda I _N ) {\mathbf W} - {\rm Cov} (y (t) , \mathbf{x} (t)) = 0,\\
a + {\mathbf W} ^\top  {\rm E} \left[ \mathbf{x} (t)\right] - {\rm E}\left[ y (t) \right] =0.
\end{align*}
These equations  yield  the pair $({\mathbf W}_{\rm out}, a_{\rm out})$ that minimizes $g({\mathbf W},a)$. We now use the fact that $\{ \mathbf{x} (t)\}_{t \in  \mathbb{Z}} $ is the unique stationary solution of VAR(1) approximating system~\eqref{VAR model reservoir} for the TDR \eqref{VAR model reservoir}  and  hence obtain 
\begin{align}
\label{linear_cov_system_gamma}
 {\mathbf W}_{\rm out} = & (\Gamma(0) + \lambda I _N ) ^{-1}{\rm Cov} (y (t) , \mathbf{x} (t)),\\
\label{linear_cov_system_c}
a_{\rm out} = & {\rm E}\left[ y (t) \right] - {\mathbf W}_{\rm out} ^\top \boldsymbol{\mu} _x,  
\end{align}
where $\boldsymbol{\mu} _x $ is provided in \eqref{mu_x}, $\Gamma(0) \in \mathbb{S} _{ N}$ is determined by the generalized Yule-Walker equations in \eqref{Gamma0} and ${\rm Cov}(y (t), \mathbf{x} (t))$ is a vector in $\mathbb{R} ^N $ that has to be determined for every specific memory task $H$. 
Additionally, it is easy to verify that the error committed by the reservoir when using the optimal readout is 
\begin{align}
\label{error rc optimal}
&{\rm E}\left[\left( {\mathbf W}_{{\rm out}} ^{\top} \cdot {\bf x}(t)  + a_{{\rm out}} -{ y}(t)\right)^2\right]= {\mathbf W}_{{\rm out}}^\top \Gamma (0)  {\mathbf W}_{{\rm out}}+{\rm var} \left(y (t)\right) - 2{\mathbf W}_{{\rm out}}^\top  {\rm Cov} (y (t) , \mathbf{x} (t))  \nonumber \\
&={\rm var} \left(y (t)\right) -{\mathbf W}_{{\rm out}}^\top (\Gamma (0) + 2\lambda \mathbb{I} _N) {\mathbf W}_{{\rm out}} \nonumber \\
	&={\rm var} \left(y (t)\right)-{\rm Cov} (y (t) , \mathbf{x} (t))^\top(\Gamma(0) + \lambda \mathbb{I} _N ) ^{-1}(\Gamma (0) + 2\lambda \mathbb{I} _N) (\Gamma(0) + \lambda \mathbb{I} _N ) ^{-1}{\rm Cov} (y (t) , \mathbf{x} (t)).
\end{align}
The $H$-memory capacity $C_H( \boldsymbol{\theta}, {\bf c}, \lambda)$ of a reservoir computer constructed using a nonlinear kernel $f$ with parameters $\boldsymbol{\theta}$, an input mask ${\bf c}  $, and regularizing ridge parameter $\lambda$, is defined as one minus the normalized mean square error committed at the time of accomplishing the memory task $H$. Expression~\eqref{error rc optimal} shows that when the RC is approximated by the VAR(1) model~\eqref{VAR model reservoir}, the corresponding {\bf $H$-memory capacity} can be approximated by
\begin{equation}
\label{capacity formula_suppl}
C_H( \boldsymbol{\theta}, {\bf c}, \lambda)=\frac{{\rm Cov} (y (t) , \mathbf{x} (t))^\top(\Gamma(0) + \lambda \mathbb{I} _N ) ^{-1}(\Gamma (0) + 2\lambda \mathbb{I} _N) (\Gamma(0) + \lambda \mathbb{I} _N ) ^{-1}{\rm Cov} (y (t) , \mathbf{x} (t))}{{\rm var} \left(y (t)\right)}
\end{equation}
Since the normalized error coming from the expression~\eqref{error rc optimal} is clearly bounded between zero and one, it is also clear that:
\begin{equation*}
0 \leq C_H( \boldsymbol{\theta}, {\bf c}, \lambda)\leq 1.
\end{equation*}
We emphasize that in order to evaluate~\eqref{capacity formula_suppl} for a specific memory task, only ${\rm Cov} (y (t) , \mathbf{x} (t)) $ and ${\rm var} \left(y (t)\right) $ need to be computed since the autocovariance $\Gamma (0)  $ is fully determined by~\eqref{Gamma0} once the reservoir and the equilibrium $\mathbf{x}_0 $ around which we operate have been chosen. 

Once a specific reservoir and task $H$ have been fixed, the capacity function $C_{H}(\boldsymbol{\theta}, {\bf c}, \lambda )$ can be explicitly written down and it can hence be used to find reservoir parameters ${\boldsymbol{\theta}  }_{{\rm opt}} $ and an input mask ${\bf c} _{\rm opt}$ that maximize it,  by solving the optimization problem
\begin{equation}
\label{optimal theta and c}
({\boldsymbol{\theta}  }_{{\rm opt}}, {\bf c} _{\rm opt}):=\mathop{\rm arg\, max}_{{\boldsymbol{\theta} }   \in \mathbb{R} ^K,  {\bf c} \in\mathbb{R}^N   } C_{H}(\boldsymbol{\theta}, {\bf c}, \lambda ).
\end{equation}

\medskip

\noindent {\bf Two specific memory tasks.} In the following paragraphs we spell out the computation of ${\rm Cov} (y (t) , \mathbf{x} (t)) $ and ${\rm var} \left(y (t)\right) $ necessary to evaluate the memory capacity formula~\eqref{capacity formula_suppl} for the two most basic memory tasks, namely, the linear and the quadratic ones.

\begin{description}

\item [(i)] {\bf The $h$-lag linear memory task.} The  linear memory task is determined by the  linear task functions $H: \mathbb{R}^{h+1}\rightarrow \mathbb{R}$ that we now describe. First, let  ${\bf z}^h(t):= \left(z (t), z (t-1), \ldots,z(t-h)\right)^\top $ and let $\mathbf{L}  \in \mathbb{R}^{h+1}$. We then set $H({\bf z}^h(t)) := \mathbf{L} ^\top {\bf z}^h(t) $. In order to evaluate the $h$-lag memory capacity using formula \eqref{capacity formula_suppl}, we need to evaluate ${\rm var} \left(y (t)\right) $ and  ${\rm Cov} (y (t) , \mathbf{x} (t)) $ with $y(t) := H({\bf z}^h(t))$. 

First, since $ \left\{ z(t)\right\} _{t \in  \mathbb{Z} }\sim {\rm IID} (0, \sigma _z ^2 )$,  we then immediately obtain that
\begin{equation}
\label{var y linear}
{\rm var} ({y} (t) ) = \sigma _{{z} } ^2 ||\mathbf{L}  || ^2.
\end{equation} 
Next, we use the so called MA($\infty$)-representation of the VAR(1) in  \eqref{VAR model reservoir}, namely, 
\begin{equation}
\label{MA_infty}
(\mathbf{x} (t) - \boldsymbol{\mu}_{ {x} }) = \sum^{\infty}_{i = 0} \Psi _i \boldsymbol{\rho} ({t-i}), 
\end{equation}
with $ \boldsymbol{\mu}_{{ x }}$ as in \eqref{mu_x},  $\boldsymbol{\rho}  (t) := \boldsymbol{\varepsilon} (t) - \boldsymbol{\mu}_{{ \varepsilon }}$, $\boldsymbol{\mu}_{{ \varepsilon }}$ defined  in \eqref{mu_eps_res_main_suppl},  $\Psi_i =A(\mathbf{x} _0,  \boldsymbol{\theta} )^i$, and  $A(\mathbf{x} _0,  \boldsymbol{\theta} )$ the connectivity matrix of the discretized nonlinear TDR provided in \eqref{DxF}. Using \eqref{MA_infty}, we compute
\begin{align}
&{\rm Cov} ({y} (t), x_i (t)) = {\rm Cov} (\mathbf{L}   ^\top \mathbf{z}^h (t), x _i (t)) =\sum^{h+1}_{j = 1} L _j {\rm Cov} (z(t-j+1), x_i(t)) \nonumber \\
&= \sum^{h+1}_{j = 1} \sum^{\infty}_{k = 0} \sum^{N}_{r = 1}L _j (A(\mathbf{x} _0,  \boldsymbol{\theta} )^k )_{ir} {\rm E} \left[ z({t-j+1}) \rho_r ({t-k}) \right]    \nonumber \\&= \sum^{h+1}_{j = 1}  \sum^{N}_{r = 1} L_j (A(\mathbf{x} _0,  \boldsymbol{\theta} ) ^{j-1} )_{ir} {\rm E} \left[ z({t})  (\varepsilon _r({t})- z(t) (\boldsymbol{\mu}_{ \varepsilon }) _r \right] , \enspace {\rm with} \enspace i \in \{1, \ldots, N\}, \label{cov_computation}
\end{align} 
which immediately yields that
\begin{equation}
\label{cov xy linear} 
{\rm Cov} ({y} (t), x_i (t)) = (1 - { e}^{- \xi }) \sum^{h+1}_{j = 1}  \sum^{N}_{r = 1} L _j (A(\mathbf{x} _0,  \boldsymbol{\theta} )^{j-1} )_{ir} p_{R} (\mu_{z}, c _1 , \dots, c _r),\enspace {\rm with} \enspace i \in \{1, \ldots, N\},
\end{equation}
where the polynomial $p_R$ on the variable $x$ is defined by $p_{R} (x, c _1 , \dots, c _r):=x \cdot q_{R} (x, c _1 , \dots, c _r)$ and its evaluation follows the same convention as in \eqref{mu_eps_res_main_suppl}.
The expressions \eqref{var y linear} and \eqref{cov xy linear} can be readily substituted in \eqref{capacity formula_suppl} in order to obtain an explicit expression for  capacity $C_{H}(\boldsymbol{\theta}, {\bf c}, \lambda )$ associated to the $h$-lag linear memory task as a function of the reservoir parameters $\boldsymbol{\theta} $ and the input mask ${\bf c}$. This expression can be subsequently treated as in \eqref{optimal theta and c} in order to determine optimal architecture parameters for this particular task.
\medskip

\item [(ii)] {\bf The $h$-lag quadratic memory task.}
In this case we use a quadratic task function $H: \mathbb{R}^{h+1}\rightarrow \mathbb{R}$ of the form 
\begin{equation}
\label{quadratic memory_suppl}
H({\bf z}^h(t)) := \mathbf{z} ^h (t) ^{\top} Q \mathbf{z} ^h (t) = \sum^{h+1}_{i = 1} \sum^{h+1}_{j = 1} Q _{ij} z({t-i+1}) z({t-j+1}),
\end{equation}
for some symmetric matrix $Q \in  \mathbb{S} _{h+1}$. Analogously to the linear task case, in order to evaluate the memory capacity  associated to $H$, we have to derive  explicit expressions for ${\rm var} \left(y (t)\right) $ and  ${\rm Cov} (y (t) , \mathbf{x} (t)) $ with $y(t) := H({\bf z}^h(t))$. The same computations as in the case of the linear task apply. First, if $ \left\{ z(t)\right\} _{t \in  \mathbb{Z} }\sim {\rm IID} (0, \sigma _z ^2 )$,  we can immediately write
\begin{equation}
\label{mu y quadratic}
{\rm E} \left[ {y} (t) \right] = \sigma _{{z} } ^2 {\rm tr}(Q),
\end{equation}
and
\begin{eqnarray}
{\rm E} \left[ {y} (t) ^2\right] &=& \sum_{i=1}^{h+1}\sum_{j=1}^{h+1}\sum_{k=1}^{h+1}\sum_{l=1}^{h+1} {\rm E}\left[Q_{ij}Q_{kl}z(t - i + 1)z(t - j + 1)z(t - k + 1)z(t - l + 1)\right]\notag\\
	&= &\sum_{i=1}^{h+1} Q_{ii}^2{\rm E}\left[z(t - i + 1)^4\right]
	+4\sum_{i=1}^{h+1}\sum_{j>i}^{h+1}Q_{ij}^2 {\rm E}\left[z(t - i + 1)^2 z(t - j + 1)^2\right]\notag\\
	& & +2\sum_{i=1}^{h+1}\sum_{j>i}^{h+1}Q_{ii} Q_{jj} {\rm E}\left[z(t - i + 1)^2 z(t - j + 1)^2\right]\notag\\
	&= &\mu _z ^4\sum_{i=1}^{h+1} Q_{ii}^2
	+4\sigma_z ^4\sum_{i=1}^{h+1}\sum_{j>i}^{h+1}Q_{ij}^2 
	+2\sigma_z ^4\sum_{i=1}^{h+1}\sum_{j>i}^{h+1}Q_{ii} Q_{jj}.\label{expectation square}\\
\end{eqnarray}
Analogously, by~\eqref{mu y quadratic},
\begin{equation}
\label{expectation square 1}
{\rm E} \left[ {y} (t) \right]^2=\sigma _{{z} } ^4 {\rm tr}(Q)^2=\sigma _{{z} } ^4 \left(\sum_{i=1}^{h+1} Q_{ii}^2+2  \sum_{i=1}^{h+1}\sum_{j>i}^{h+1}Q_{ii} Q_{jj} \right).
\end{equation}
Hence, if we put together~\eqref{expectation square} and~\eqref{expectation square 1}, we obtain
\begin{equation}
\label{var y quadratic}
{\rm var} ({y} (t) ) = {\rm E} \left[ {y} (t)^2 \right] -{\rm E} \left[ {y} (t) \right]^2 
  = (\mu _z ^4- \sigma _z ^4)\sum_{i=1}^{h+1} Q_{ii}^2+4\sigma_z ^4\sum_{i=1}^{h+1}\sum_{j>i}^{h+1}Q_{ij}^2. 
\end{equation}
Recall that for Gaussian variables, that is  $ \left\{ z(t)\right\} _{t \in  \mathbb{Z} }\sim{\rm IN} (0, \sigma _z ^2 )$, we have that $\mu _z ^4= 3 \sigma _z ^2$, and hence in that case 
\begin{equation}
\label{var y quadratic gaussian}
{\rm var} ({y} (t) ) 
  = 2\sigma _z ^4\left(\sum_{i=1}^{h+1} Q_{ii}^2+2\sum_{i=1}^{h+1}\sum_{j>i}^{h+1}Q_{ij}^2\right)=2\sigma _z ^4\sum_{i=1}^{h+1}\sum_{j=1}^{h+1}Q_{ij}^2. 
\end{equation}

Regarding the computation of the covariance and analogously to the case of the linear $h$-lag memory task, we use the MA($\infty$) representation of the VAR(1) model of the TDR in  \eqref{VAR model reservoir} and write
\begin{align*}
\label{cov_computation}
{\rm Cov} ({y} (t), x_i (t)) = &\sum^{h+1}_{j = 1} \sum^{h+1}_{k = 1}\sum^{\infty}_{l = 1}Q_{jk} {\rm E} \left[ z({t-j+1}) z({t-k+1}) ( A(\mathbf{x} _0,  \boldsymbol{\theta} ) ^l \boldsymbol{\rho} ({t-l})) _i \right] \nonumber\\
=&  \sum^{\infty}_{l = 0} \sum^{h+1}_{j = 1} \sum^{h+1}_{k = 1}\sum^{N}_{ r = 1} Q _{jk} (A(\mathbf{x} _0,  \boldsymbol{\theta} )^l)_{ir} \Big\{ {\rm E} \left[ z({t-j+1}) z({t-k+1}) \varepsilon_r  ({t-l})  \right] \\
&- (\boldsymbol{\mu} _{{ \varepsilon }})_r {\rm E} \left[  z({t-j+1}) z({t-k+1})\right]  \Big\}, \enspace {\rm with} \enspace i \in \{1, \ldots, N\},\nonumber
\end{align*} 
which leads to the following result:
\begin{equation}
\label{cov xy quadratic}
{\rm Cov} ({y} (t), x_i (t)) =   (1 - { e}^{- \xi }) \sum^{h+1}_{j = 1}  \sum^{N}_{r = 1} Q_{jj} (A(\mathbf{x} _0,  \boldsymbol{\theta} ) ^{j-1} )_{ir}(s_{R} (\mu_{{z} }, c _1 , \dots, c _r) - \sigma _{z}^2q_{R}  (\mu_z, c _1 , \dots, c _r)),
\end{equation} 
 { with} $ i \in \{1, \ldots, N\}$.  In this relation the polynomial $s_R$ on $x$ is defined as
$s_{R} (x, c _1 , \dots, c _r):=x ^2  \cdot q_{R} (x, c _1 , \dots, c _r)$ and is evaluated following the same convention as in \eqref{mu_eps_res_main_suppl} but taking $x^2$ instead of $x$.

Again, we conclude by noticing that the expressions \eqref{cov xy quadratic} and \eqref{var y quadratic} substituted in \eqref{capacity formula_suppl} provide an explicit formula for  capacity $C_{H}(\boldsymbol{\theta}, {\bf c}, \lambda )$ associated to the $h$-lag quadratic memory task as a function of the reservoir parameters $\boldsymbol{\theta} $ and the input mask ${\bf c}$ and hence it can be readily used to solve the optimization problem in \eqref{optimal theta and c}.

An observation that is worth to be pointed out is that only the diagonal elements in $Q$ intervene in the covariance~\eqref{cov xy quadratic} while all its entries are present in the variance~\eqref{var y quadratic}. When these two quantities are substituted in the memory capacity formula \eqref{capacity formula_suppl} it can be seen that by choosing sufficiently high off-diagonal entries in $Q$, the capacity of the reservoir can be made arbitrarily small which shows a structural limitation of the architecture that we are considering that can only be fixed by using alternative signal feeding schemes.
\end{description}

 \section{Equilibria of the continuous and the discrete time models for the TDR and their stability}
\label{Equilibria of time-delay differential equation and those of its discretization}
As we already explained, the linearization of the reservoir map at a stable fixed point is at the core of the developments in this paper. That is why in this section we carry out a detailed study of the stability properties of the equilibria of the time-delay differential equation  \eqref{time delay equation} and of the fixed points of its corresponding discrete-time approximation  \eqref{vector discretized reservoir}.
More specifically, we provide sufficient stability conditions and we show that our results exhibit a remarkable consistence regardless of the use of the continuous or of the discrete time schemes.

\subsection{Stationary solutions of time-delay differential equations and their stability}
\label{Solutions of time-delay differential equations and their stability} 
We  start  by recalling some basic facts about the properties of the solutions of the time-delay differential equations and their stability. Let $\tau \in \mathbb{R} ^+$ be a fixed delay and consider a {\bf time-delay map} 
\begin{align}
\label{map X general}
\begin{array}{ccccccc}X:&{C} ^1([ - \tau , 0], \mathbb{R})  \times \mathbb{R}&\longrightarrow&\mathbb{R}\\&({\gamma  }, t)&\longmapsto&X( \gamma , t).\end{array}
\end{align}
Additionally, for any $t \in \mathbb{R} $ define the {\bf shift operator} 
\begin{align}
\label{shift operator}
\begin{array}{ccccccc}S_t:&{C} ^1([- \tau +t, t], \mathbb{R})   &\longrightarrow&{C} ^1([ - \tau , 0], \mathbb{R})\\&\gamma&\longmapsto&\gamma \circ \lambda_{t},\end{array}
\end{align}
where $\lambda_{t}$ is the translation operator by $t \in \mathbb{R} $, that is, $\lambda_{t}(s):=s+t$, for any $s \in \mathbb{R}$. Let now $\gamma \in {C} ^1([ - \tau , +\infty), \mathbb{R}) $ be a differentiable curve. We say that $\gamma$ is a solution of the {\bf time-delay differential equation (TDDE) determined by $X$} when the equality 
\begin{equation}
\label{global notation for tdde}
\dot \gamma (t)= X(S _t \circ \gamma|_{[-\tau+t, t]}, t)
\end{equation} 
holds for any $t \in [0, +\infty) $. Note that the TDDE~\eqref{time delay equation} that is at the core of this paper, namely 
\begin{equation}
\label{time delay equation 2 recall_suppl}
\dot{x}(t) = -x(t) + f(x(t - \tau ), I(t), \boldsymbol{\theta} ),
\end{equation}
can be encoded as in~\eqref{global notation for tdde} by using the time-delay map $X$ given by
\begin{align}
\label{map X nonlinear TDR}
\begin{array}{ccccccc}X:&{C} ^1([ - \tau , 0], \mathbb{R}) \times \mathbb{R}&\longrightarrow&\mathbb{R}\\&({\gamma  }, t)&\longmapsto&- \gamma (0) +f(\gamma (- \tau ), I(t), \boldsymbol{\theta} ).
\end{array}
\end{align}

\begin{definition}
We say that the time-delay map $X$ is locally Lipschitzian on the open set $\Omega \subset {C} ^1([ - \tau , 0], \mathbb{R})  \times \mathbb{R}$ if it is Lipschitzian in any compact subset of $\Omega  $, that is, for any compact subset $\Omega _0 $ of $\Omega  $ there exists a constant $K \in \mathbb{R} ^+$ such that for all $(\gamma _1,t) $ and $(\gamma _2 , t)$ in $\Omega _0$ one has   
\begin{equation}
\label{Lipschitz}
|X(\gamma _1 , t) - X(\gamma _2 , t)| < K ||\gamma _1 -\gamma _2 ||_{\infty}.
\end{equation}
\end{definition}
\begin{theorem}[Existence and uniqueness of solutions]
Let $X$ be a continuous and locally Lipschitzian time-delay map in  ${C} ^1([ - \tau , 0], \mathbb{R})  \times \mathbb{R}$. Then, for any $\phi \in {C} ^1([ - \tau , 0], \mathbb{R})$ there exists a unique $\Gamma _{\phi } \in {C}^1 ([- \tau , +\infty), \mathbb{R} )$ such that  
\begin{align}
\label{TDDE} 
&\Bigg\{
\begin{array}{rll} 
\Gamma _{\phi } (t)& = \phi (t), \enspace &{\rm for} \enspace {\rm any } \enspace t \in [- \tau , 0]\\
\dot{ \Gamma} _{\phi } (t) &= X(S_t \circ \Gamma _{\phi }|_{[- \tau+t,t]}  , t), \enspace &{\rm for } \enspace{\rm any} \enspace  t\in (0,+\infty].
\end{array}
\end{align}
We say that $\Gamma _{\phi}$ is the {\bf solution} of the time-delay differential equation determined by $X$ with initial condition $\phi$, or simply the solution through $\phi$.	The associated {\bf flow} is defined as the map
\begin{align}
\label{flow}
\begin{array}{cccc}
F: & [- \tau , +\infty) \times {C} ^1([ - \tau , 0], \mathbb{R}) &\longrightarrow &\mathbb{R} \\
&(t, \phi )&\longmapsto& \Gamma_{\phi } (t)
\end{array}
\end{align}
and note that $F_{\boldsymbol{\cdot}} ( \phi ) \in {C}^1([- \tau , +\infty),\mathbb{R} )$.
\end{theorem}

We now recall also some basic notions of  stability of common use in the TDDE context; see \cite{Hale:functional_differential_equations} and  \cite{she:book} for details. Let $x _0 \in \mathbb{R} $ and let $\phi_{x _0} \in   {C} ^1([ - \tau , 0], \mathbb{R}) $ be the constant curve at $x _0 $. We say that the point  $x _0 $ is an {\bf equilibrium} of the TDDE determined by the time-delay map and with flow $F $ whenever $F _t(\phi_{x _0}) =x _0$, for any $t \in [- \tau , +\infty)$. The equilibrium $x _0 $ is said to be {\bf stable} (respectively {\bf asymptotically stable}) if for any $\epsilon >0  $ there exists a $\delta(\epsilon)> 0 $ such that for any $\phi \in   {C} ^1([ - \tau , 0], \mathbb{R})$  with $\| \phi- \phi_{x _0}\| _{\infty}<\delta(\epsilon) $, we have that $|F _t(\phi)- x _0|< \epsilon $, for any $t \in [- \tau , +\infty) $ (respectively $\lim\limits _{t \rightarrow \infty}F _t(\phi)= x _0$). 

The following stability criterion is an extension of Lyapunov's Second Method to the TDDE context due to Krasovskiy~\cite{krasovskiy:book}. We state it using our notation since it will be used in the sequel.

\begin{theorem}[Lyapunov-Krasovskiy stability theorem]
\label{Lyapunov-Krasovskiy stability theorem} 
let $x _0 \in \mathbb{R} $ be an equilibrium of the time-delay differential equation 
\eqref{global notation for tdde} with flow $F:   [- \tau , +\infty)\times {C} ^1([ - \tau , 0], \mathbb{R})) \longrightarrow \mathbb{R}$. Let  $u$, $v$, $w: \overline{{\mathbb{R}} ^+}   \longrightarrow \overline{\mathbb{R} ^+}$ be continuous nondecreasing functions such that  $u(0)=v(0) = 0$ and $u(t), v(t), w(t)>0$ for any $t \in (0, + \infty)$. If there exists a continuously differentiable functional $V$ 
\begin{align}
\label{V Lyapunov-Krasovskiy functional}
\begin{array}{ccccccc}V:& {C}^1 ([ - \tau , +\infty), \mathbb{R} )  \times \mathbb{R} &\longrightarrow\mathbb{R} \end{array}
\end{align}
such that  for any $\phi \in {C} ^1([ - \tau , 0], \mathbb{R})) $ and any $t \in [ 0, +\infty)$ satisfies that 
\begin{description}
\item[(i)] $u(|\phi (0)|) \le V(F_{\boldsymbol{\cdot}}( \phi ), t) \le v(||\phi ||_\infty)$,
\item [(ii)] $\dot{V}(F_{\boldsymbol{\cdot}}( \phi ), t) := \dfrac{d}{dt} V(F_{\boldsymbol{\cdot}}( \phi ),t) \leq -w(| \phi (0)|)$, 
\end{description}
then $x_0$ is asymptotically stable. If $w(t) \ge 0$ then $x_0$ is just  stable. A functional $V$ that satisfies these conditions is called a {\bf Lyapunov-Krasovskiy functional}.
\end{theorem}

\subsection{Equilibria of the  reservoir time-delay  equation and their stability}

We now use Theorem~\ref{Lyapunov-Krasovskiy stability theorem} to establish sufficient conditions for the stability of the equilibria  of the TDDE \eqref{time delay equation} at the core of the paper, namely, 
\begin{equation}
\label{time delay equation 2 recall}
\dot{x}(t) = x(t) + f(x(t - \tau ), I(t), \boldsymbol{\theta} ).
\end{equation}
where $f$ is the nonlinear kernel of the TDR. The main tool in the application of that result is the use of a Lyapunov-Krasovskiy functional of the form
\begin{align}
\label{V for nonlinear TDR}
\begin{array}{ccccccc}V:&{C} ^{1}([ - \tau , +\infty],\mathbb{R} )  \times \mathbb{R}&\longrightarrow&\mathbb{R}\\&({x_{\phi }  }, t)&\longmapsto&\dfrac{1}{2}x_{\phi } (t)^2 +{ \mu } \int ^t_{t - \tau } x_{\phi } (s)^2 ds,\end{array}
\end{align}
where $\mu \in \mathbb{R} ^+$ and $x_{\phi} = F_{\boldsymbol{\cdot}}(\phi )$ for some initial curve $\phi \in {C} ^1([ - \tau , 0], \mathbb{R})$. See~\cite{krasovskiy:book}, \cite{Hale:functional_differential_equations} and \cite{she:book} for the extensive discussion.

\begin{theorem}
\label{theorem stability continuous} Let $ x_0 $ be an equilibrium of the time-delay differential equation \eqref{time delay equation 2 recall}
in autonomous regime, that is, when $I(t) = 0$, and suppose that there exists $\varepsilon >0$ and $k_{\varepsilon } \in \mathbb{R} $ such that one of the following conditions holds
\begin{description}
\item[(i)] $f(x + x _0, 0, \boldsymbol{\theta}  )\le k_{\varepsilon }x + x _0 $ for all $x \in \left( - \varepsilon , \varepsilon\right) $
\item[(ii)] $\dfrac{f(x + x _0, 0,\boldsymbol{\theta}  ) - x_0 }{x} \le k_{\varepsilon }$ for all $x \in \left( - \varepsilon , \varepsilon\right) $.
\end{description} 
If $|k_{\varepsilon }|<1$ then $x_0 $ is {\bf  asymptotically stable}. If $|k_{ \varepsilon }| \le  1$ then $x _0 $ is {\bf  stable}.
\end{theorem}
{\bf{Proof.}} Notice first that the equilibria $x _0 $ in the statement are characterized by the equality $f(x _0, 0, \boldsymbol{\theta}  )= x _0 $.
Consider now the Lyapunov-Krasovskiy functional introduced in \eqref{V for nonlinear TDR}.
It is easy to see that since $V(x_{\phi}, t)$ is positive it satisfies  condition {\bf (i)} in Theorem~\ref{V Lyapunov-Krasovskiy functional}. We will now show that any of the two conditions in the statement imply that condition {\bf (ii)}  in Theorem~\ref{V Lyapunov-Krasovskiy functional} are satisfied and hence guarantee the stability of $x _0$. We start by writing
\begin{align}
\label{dVdt}
\dfrac{d}{dt}V(x_{\phi},t) &= x_{\phi}(t) \dot{x}_{\phi}(t) + \mu \left( x_{\phi}(t) ^2 - x_{\phi}(t - \tau ) ^2 \right) \nonumber \\
&= -x_{\phi}(t) ^2 + x_{\phi}(t) f(x_{\phi}(t - \tau ), 0, \boldsymbol{\theta} ) +\mu\left( x_{\phi}(t) ^2 - x_{\phi}(t- \tau ) ^2  \right).
\end{align}
We now distinguish two cases, namely, when $x _0= 0 $ and when $x _0\neq 0 $.

\medskip

\noindent {\bf Case $x _0 =0 $.}
Suppose that $x_0 =0$ is a solution of the TDDE \eqref{time delay equation 2 recall}. Under the hypothesis {\bf (i)} in the statement, in the case of the trivial solution  $x_0=0$ there exists $\varepsilon >0$ and $k_{\varepsilon }>0$ such that $f(x, 0, \boldsymbol{\theta} ) \le  k_{\varepsilon }x$ for all $x \in (- \varepsilon , \varepsilon )$, and hence from~\eqref{dVdt} we can conclude that 
\begin{align}
\label{dVdtQ}
\dfrac{d}{dt}V(x_{\phi},t) & \le -x_{\phi}(t) ^2 + k_{\varepsilon } x_{\phi}(t) x_{\phi}(t - \tau ) +\mu\left( x_{\phi}(t) ^2 - x_{\phi}(t- \tau ) ^2  \right)  \nonumber \\
&= \left( x_{\phi}(t), x_{\phi}(t - \tau )\right) Q\left( x_{\phi}(t), x_{\phi}(t - \tau )\right) ^\top 
\end{align}
 with $$Q := \left( \begin{array}{cc} \mu -1 & {k_{\varepsilon }}/{2}\\{k_{\varepsilon }}/{2}& -\mu \end{array}\right). $$
Expression \eqref{dVdtQ} is negative for any $\left( x_{\phi}(t), x_{\phi}(t - \tau )\right) $  if the matrix $Q$ is negative definite which by the Sylvester's law amounts to $\mu < 1$ and $k_{\varepsilon } ^2 < -4\mu ( \mu - 1)$. Since $-4\mu (\mu -1)$ has a maximum at $\mu ={1}/{2}$ for which $-4 \mu (\mu -1) = 1$, we obtain from Theorem~\ref{Lyapunov-Krasovskiy stability theorem} that the optimal  sufficient condition for asymptotic stability of $x_0$  is $|k_{\varepsilon} |  < 1$ as required. Analogously, by Theorem~\ref{Lyapunov-Krasovskiy stability theorem}, a sufficient condition for $x_0=0$ to be stable is the non-positivity of expression  \eqref{dVdtQ} or, equivalently, the negative semi-definiteness of $Q$ which amounts to $\mu \le 1$ and $k_{\varepsilon } ^2 \le -4\mu ( \mu - 1)$. Hence the optimal  sufficient condition for the stability of $x_0=0$ is $|k_{\varepsilon} |  \le  1$ as required.
 
 Consider   \eqref{dVdt} again, now under the hypothesis {\bf (ii)} of the statement of the theorem. In the case of the trivial solution it implies that there exists $\varepsilon >0$ and $k_{\varepsilon }>0$, such that $\dfrac{f(x, 0, \boldsymbol{\theta} )}{x} \le  k_{\varepsilon }$ for all $x \in (- \varepsilon , \varepsilon )$. In order to ensure the asymptotic stability of the trivial solution using  Theorem~\ref{Lyapunov-Krasovskiy stability theorem}, we need to find conditions under which the expression \eqref{dVdt} is negative, that is $\dfrac{d}{dt}V(x_{\phi },t)< 0$. We proceed by first multiplying both sides of this   inequality by the positive quantity $\dfrac{1}{x_{\phi }(t - \tau) ^2}$. We obtain
\begin{equation*}
-\dfrac{x_{\phi }(t) ^2}{x_{\phi }(t - \tau) ^2} + \dfrac{x_{\phi }(t) f(x_{\phi } (t- \tau ), 0, \boldsymbol{\theta} )}{x_{\phi }(t - \tau ) ^2 } +\mu\left( \dfrac{x_{\phi }(t) ^2}{x_{\phi }(t - \tau) ^2 } - 1  \right)  < 0.
\end{equation*}
Then due to the hypothesis {\bf (ii)} of the theorem, a sufficient condition for this inequality to hold is 
\begin{equation}
\label{ii inequality} 
-\dfrac{x_{\phi }(t) ^2}{x_{\phi }(t - \tau) ^2} + k_{\varepsilon } \dfrac{x_{\phi }(t)}{x_{\phi }(t - \tau )} +\mu\left( \dfrac{x_{\phi }(t) ^2}{x_{\phi }(t - \tau) ^2 } - 1  \right)  < 0.
\end{equation}
Notice that when $x_{\phi }(t) = x_{\phi }(t - \tau ) = x _0 $, this inequality is always satisfied provided that $k_{\varepsilon }\le 1$. Hence in order for \eqref{ii inequality} to hold, it suffices that the polynomial on $z$
\begin{equation}
\label{polynomial for ii}
-z ^2 + k_{\varepsilon } z + \mu ( z ^2 - 1) = (\mu -1) z ^2 + k_{\varepsilon }z - \mu
\end{equation}
has no real roots, which happens, as in point {\bf (i)} of the statement when $ k_{ \varepsilon } ^2 < -4 \mu ( \mu -1 )$. Proceeding analogously as under assumption {\bf(i)} of the theorem, we obtain $|k_{\varepsilon }|<1$ (respectively $|k_{\varepsilon }|\le 1$) as the sufficient condition for the asymptotic stability (respectively stability) of $x_0=0$, as required. 

\medskip

\noindent {\bf Case $x _0\neq 0 $.}
Suppose now that $x_0 \ne 0$ and define the new variable $y(t):=x(t) - x_0$. With this change of variables the equation  \eqref{time delay equation 2 recall} becomes 
\begin{equation*}
\label{TDE change variables}
\dot{y}(t) = \dot{x}(t) = -x_0 - y(t) + f(y( t- \tau ), 0,\boldsymbol{\theta} )
\end{equation*}
or, equivalently,
\begin{equation}
\label{TDE change variables g}
\dot{y}(t) =  - y(t) + g(y( t- \tau ), 0,\boldsymbol{\theta} ),
\end{equation}
where the function $g$ is defined as 
$g(y(t), 0,\boldsymbol{\theta} ) := f(y(t)+x_0, 0, \boldsymbol{\theta} ) - x_0$. The equation \eqref{TDE change variables} has an equilibrium at $y_0 = 0$ whose stability can be easily studied by mimicking the case $x_0 = 0$  discussed above. More specifically,  it can be shown following the same arguments that in this case the hypothesis {\bf (i)} of the statement of the theorem can be written as
\begin{equation}
\label{g hypothesis (i)}
g(y, 0, \boldsymbol{\theta} )\le k_{\varepsilon } y
\end{equation}
and $y_0 = 0$ is stable or asymptotically stable whenever $|k_{\varepsilon }| \le 1$ or $|k_\varepsilon | < 1$, respectively. The inequality \eqref{g hypothesis (i)} is equivalent to
\begin{equation*}
\label{}
f(y + x_0,0,\boldsymbol{\theta} ) - x_0 \le k_{\varepsilon } y
\end{equation*}
or
\begin{equation*}
\label{}
f(y + x_0,0,\boldsymbol{\theta} )  \le k_{\varepsilon } y + x_0,
\end{equation*}
which guarantees that a non-trivial equilibrium $x_0$ of \eqref{time delay equation 2 recall} is stable or asymptotically stable when the same conditions on $k_{\varepsilon }$ as in the trivial case are satisfied.

Finally, the hypothesis {\bf (ii)} of the statement of the theorem in the case of  \eqref{TDE change variables g} has the form 
\begin{equation}
\label{g hypothesis (ii)}
\dfrac{g(y, 0, \boldsymbol{\theta} )}{y}\le k_{\varepsilon }
\end{equation}
and $y_0 = 0$ is stable or asymptotically stable when $|k_{\varepsilon }| \le 1$ or $|k_\varepsilon | < 1$, respectively. It is easy to verify that the inequality \eqref{g hypothesis (ii)} is equivalent to
\begin{equation}
\label{}
\dfrac{f(y + x_0,0,\boldsymbol{\theta} ) - x_0}{y} \le k_{\varepsilon },
\end{equation}
which provides the same corresponding sufficient conditions on $k_{\varepsilon }$ for stability or asymptotic stability of a non-trivial equilibrium $x_0$ of \eqref{time delay equation 2 recall}, as required.
$\square$

\begin{corollary}
\label{Corollary 2} 
Let $x_0$ be an equilibrium of the TDDE \eqref{time delay equation 2 recall} and suppose that 
the  nonlinear reservoir kernel function $f$ is continuously  differentiable at $x_0$.  If  $|\partial_{x}f(x_0, 0, \boldsymbol{\theta} )| < 1$ (respectively, $|\partial_{x}f(x_0, 0, \boldsymbol{\theta} )| \le  1$), then $x_0$ is asymptotically stable (respectively, stable). 
\end{corollary}

\noindent{\bf Proof.} First, define the function
\begin{subnumcases}{g_{\varepsilon}(h):=}
\dfrac{f(x_0 +h,0, \boldsymbol{\theta} )-f(x_0,0, \boldsymbol{\theta} )}{h},& $h\neq 0, h\in (- \varepsilon , \varepsilon )$\label{corollary 2 function g}\\
\partial_{x}f(x_0, 0, \boldsymbol{\theta} ) = \lim_{h \rightarrow 0} \dfrac{f(x_0 +h,0, \boldsymbol{\theta} )-f(x_0,0, \boldsymbol{\theta} )}{h},&$h=0$.\label{corollary 2 function g1}
\end{subnumcases}
By construction, the function $g_{\varepsilon}$ is continuous in $(-\varepsilon , \varepsilon )$, that is $g_{\varepsilon} \in C^0((- \varepsilon , \varepsilon ), \mathbb{R} )$. Hence, by the Weierstrass extreme value theorem, this function reaches a maximum $k_{\varepsilon }$ in the interval $[- {\varepsilon}/{2} , {\varepsilon}/{2} ]$, that is,
\begin{equation}
\label{corollary 2 function g k}
g_{\varepsilon}(h) = \dfrac{f(x_0 +h,0, \boldsymbol{\theta} )-f(x_0,0, \boldsymbol{\theta} )}{h} \le k_{\varepsilon } \quad \mbox{for any} \quad h \in [- {\varepsilon}/{2} , {\varepsilon}/{2} ].
\end{equation}
Since $x_0$ is an equilibrium, then $f(x_0,0, \boldsymbol{\theta} )  =x_0$ and the condition \eqref{corollary 2 function g k} coincides with the hypothesis {\bf(ii)} of Theorem~\ref{theorem stability continuous}. The equilibrium $x_0$ can be hence proved to be asymptotically stable (respectively, stable) if $|k_ \varepsilon|<1 $ (respectively, $|k_ \varepsilon| \le 1 $). Additionally, using~\eqref{corollary 2 function g}-\eqref{corollary 2 function g1} and the continuity of $g _ \varepsilon  $, it is easy to see that 
\begin{equation*}
\lim\limits_{\varepsilon  \rightarrow 0}k _ \varepsilon= \lim_{\varepsilon  \rightarrow 0
} g_ \varepsilon(h) = \partial_{x}f(x_0, 0, \boldsymbol{\theta} ) 
\end{equation*}
and hence the asymptotic stability (respectively, stability) of $x_0$ is guaranteed if   $|\partial_{x}f(x_0, 0, \boldsymbol{\theta} |<1 $ (respectively, $|\partial_{x}f(x_0, 0, \boldsymbol{\theta} | \le 1 $), as required.
$\square$

\medskip

We now study the equilibria and the parameter values that ensure their stability when Corollary~\ref{Corollary 2} is applied to the two nonlinear kernels that are most used in our work, that is, the Mackey-Glass~\cite{mackey-glass:paper} and the Ikeda~\cite{Ikeda1979} parametric families. 
We recall that the {\bf Mackey-Glass nonlinear kernel} is given by the  expression
\begin{equation}
\label{MG nonlinear kernel}
f(x,I, \boldsymbol{\theta})= \frac{\eta \left(x+\gamma I \right)}{1+ \left(x+\gamma I\right)^p},
\end{equation} 
where the parameter $\boldsymbol{\theta}:=(\gamma, \eta, p) \in \mathbb{R} ^3 $ is a three tuple of real values. 
The {\bf Ikeda nonlinear  kernel} corresponds to 
\begin{equation}
\label{Ikeda nonlinear kernel}
f(x,I, \boldsymbol{\theta})=\eta \sin^2(x+ \gamma I + \phi ),
\end{equation} 
where the parameter vector $\boldsymbol{\theta}:=(\gamma, \eta , \phi ) \in \mathbb{R} ^3 $. In both cases the parameter $\gamma  $ is called the {\bf input gain } and $\eta  $ the {\bf feedback gain}.
\begin{corollary}[Stability of the equilibria of the Mackey-Glass TDDE] 
\label{Stability MG continuous} 
Consider the TDDE~\eqref{time delay equation 2 recall} in the autonomous regime constructed with the Mackey-Glass kernel \eqref{MG nonlinear kernel} with $p=2$, that is,
\begin{equation}
\label{MG kernel autonomous p2}
f(x, 0, \boldsymbol{\theta} ) = \dfrac{\eta x}{1+x^2}.
\end{equation}
This TDDE exhibits two families of equilibria depending on the values of $\eta$:
\begin{description}
\item[(i)]  The trivial solution $x_0 = 0$, for any $\eta \in \mathbb{R} $. The equilibrium $x_0 = 0$ is asymptotically stable (respectively, stable) if $| \eta|<1$ (respectively, $|\eta| \le 1$).
\item[(ii)]  The non-trivial solutions $x_0 = \pm \sqrt{\eta - 1}$, for any $\eta > 1$. The equilibria $x_0 = \pm \sqrt{\eta - 1}$ are asymptotically stable (respectively, stable) whenever $1<\eta < 3$ (respectively, $1<\eta \le 3$).
\end{description} 
\end{corollary}

\noindent{\bf Proof.}
First, in order to characterize the equilibria of the time-delay differential equation \eqref{time delay equation 2 recall} with the nonlinear kernel in \eqref{MG kernel autonomous p2}, we solve $0 = -x + f(x, 0,\boldsymbol{\theta} )$ or, equivalently,
\begin{equation*}
\label{}
\dfrac{\eta x}{1+x^2} - x = 0.
\end{equation*}
A straightforward computation shows that this equality is equivalent to $x(x^2 - (\eta - 1)) = 0$ which immediately yields the two families of equilibria in the statement, namely, $x_0 = 0$, $\forall \eta \in \mathbb{R} $ and $x_0 = \pm \sqrt{\eta - 1}$ for any $\eta > 1$. We now use  Corollary~\ref{Corollary 2} of Theorem~\ref{theorem stability continuous}, in order to provide the sufficient conditions for stability and asymptotic stability of these two families. Using ~\eqref{MG  kernel autonomous p2}, we obtain that
\begin{equation}
\label{MG x derivative}
\partial_{x}f(x, 0, \boldsymbol{\theta} ) = \dfrac{\eta(1-x ^2 )}{1+ x ^2}.
\end{equation}
Then, when we evaluate this expression at the equilibria under study, we obtain:
\begin{description}
\item [(i)] for $x_0 = 0$, we have that $\partial_{x}f(x_0, 0, \boldsymbol{\theta} ) = \eta$ and hence by Corollary~\ref{Corollary 2} the trivial solution $x_0$ is asymptotically stable (respectively, stable) if $| \eta|<1$ (respectively, $|\eta| \le 1$).
\item [(ii)] for $x_0 = \pm \sqrt{\eta - 1}$ with $\eta>1$ the expression \eqref{MG x derivative} amounts to $\partial_{x}f(x_0, 0, \boldsymbol{\theta} ) = {2-\eta}$ and hence by Corollary~\ref{Corollary 2} the non-trivial solutions $x_0$ are asymptotically stable (respectively, stable) whenever $\eta \in (1,3)$ (respectively, $\eta \in (1,3]$), as required. 
\end{description}
\quad $\square$
\begin{corollary}[Stability of the equilibria of the  Ikeda TDDE] 
\label{suppl-Stability Ikeda continuous} 
Consider the TDDE~\eqref{time delay equation 2 recall} in the autonomous regime constructed with the Ikeda kernel \eqref{Ikeda nonlinear kernel}, that is,
\begin{equation}
\label{Ikeda kernel autonomous}
f(x, 0, \boldsymbol{\theta} ) = \eta \sin ^2 (x+ \phi ).
\end{equation}
The Ikeda nonlinear TDDE exhibits two families of equilibria:
\begin{description}
\item [(i)]
The trivial solution $x_0=0$ for any $\eta \in \mathbb{R} $ and $\phi = \pi n$, $n\in  \mathbb{Z} $. The equilibium $x_0 =0$ is asymptotically stable for any $\eta \in \mathbb{R} $.
\item  [(ii)]
The non-trivial equilibria $x _0$ are obtained as solutions of the equation $x_0 = \eta \sin ^2 (x_0 + \phi )$, for any $\eta \in \mathbb{R} $ and $\phi \ne \pi n$, $n\in  \mathbb{Z} $. These equilibria are asymptotically stable (respectively, stable)  whenever 
\begin{equation}
\label{stability inequality}
|\sin (2x_0 + 2\phi)| < \dfrac{1}{| \eta |} \quad (\mbox{respectively,\ } |\sin (2x_0 + 2\phi)| \le \dfrac{1}{| \eta |}). 
\end{equation}
When $| \eta | < 1$  (respectively, $| \eta |\le1$),  there exists only one non-trivial equilibrium that is always asymptotically stable (respectively, stable).
\end{description}
\end{corollary}

\noindent{\bf Proof.}  
The equilibria of the time-delay differential equation  \eqref{time delay equation 2 recall} with the Ikeda  kernel  \eqref{Ikeda kernel autonomous}, are characterized by the roots $x_0$ of the equation $0 = -x + f(x, 0,\boldsymbol{\theta} )$ or, equivalently,
\begin{equation}
\label{Ikeda equilibrium equation}
\eta \sin ^2 (x + \phi ) - x= 0.
\end{equation}
We divide the solutions of this equation into two families, namely, the trivial equilibrium $x_0 = 0$, for any $\eta \in \mathbb{R} $ and $\phi = \pi n$, $n\in  \mathbb{Z} $, and the non-trivial ones  obtained when $\phi \neq \pi n$, $n\in  \mathbb{Z} $.
Using ~\eqref{Ikeda  kernel autonomous}, we compute
\begin{equation}
\label{Ikeda x derivative}
\partial_{x}f(x, 0, \boldsymbol{\theta} ) = \eta \sin (2x +2 \phi )
\end{equation}
and evaluate it at the two families of equilibria under study.
\begin{description}
\item [(i)] For $x_0 = 0$ the expression \eqref{Ikeda x derivative} yields $\partial_{x}f(x_0, 0, \boldsymbol{\theta} ) =\eta \sin (2\phi ) \equiv 0$, since $\phi = \pi n$, $n \in  \mathbb{Z} $. Hence by Corollary~\ref{Corollary 2} the trivial solution $x_0$ is always asymptotically stable.
\item [(ii)] For non-trivial equilibria $x_0$,  the expression \eqref{Ikeda x derivative} amounts to $\partial_{x}f(x_0, 0, \boldsymbol{\theta} ) = \eta \sin (2 x_0 + 2\phi )$ and hence by Corollary~\ref{Corollary 2} the non-trivial solutions $x_0$ are asymptotically stable (respectively, stable) whenever $|\sin (2 x_0 + 2\phi )| < \dfrac{ 1}{|\eta |}$ (respectively, $|\sin (2 x_0 + 2\phi )| \le \dfrac{ 1}{|\eta |}$). We now consider the case $|\eta | < 1$ (respectively, $|\eta | \le 1$); in that situation the stability inequalities~\eqref{stability inequality} always hold true but it remains to be shown that only one equilibrium exists. That claim is a consequence of the following lemma.  

\begin{lemma}
If $|\eta | < 1$, then the equation \eqref{Ikeda equilibrium equation} has at most one root.
\end{lemma}
{\bf Proof of Lemma.} Consider the function $g(x):=\eta \sin ^2 (x + \phi ) - x$. As $g'(x)=\eta \sin (2x+2 \phi ) -1 $, we have that if $|\eta | \le 1$, then $g'(x) \le \eta - 1\le 0 $ for any $x \in \mathbb{R} $. The function $g(x)$ is hence a monotonously decreasing function and intersects the $OX$ axis in at most one point. Since $g(0) >0$ (recall that in this case $\phi \ne \pi n$, $n\in  \mathbb{Z} $) and for any $x> \eta $ we have that $g (x)<0 $, we conclude that $g(x)$ intersects the $OX$ axis in exactly one point, as required.
$\square$
\end{description}

\begin{figure}[ht]
\hspace{0.8cm}\includegraphics[scale=.65]{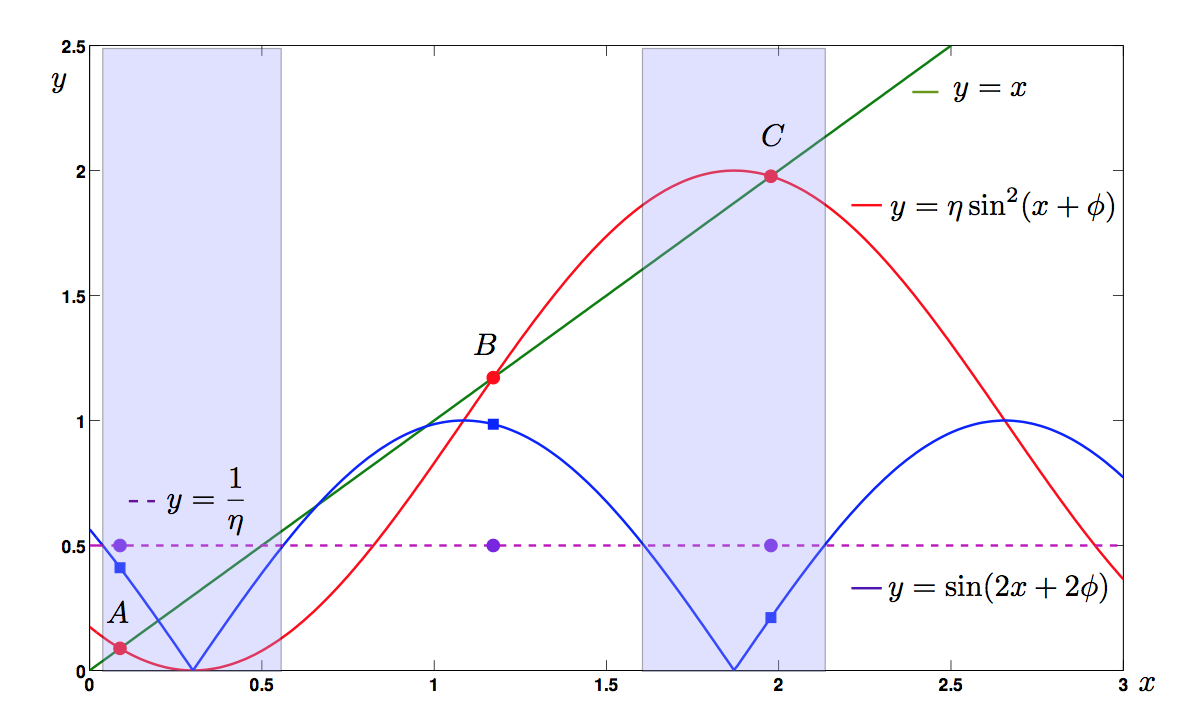}
\caption{Illustration of the statement of Corollary~\ref{suppl-Stability Ikeda continuous}. The parameter vector $\boldsymbol{\theta} $ in \eqref{Ikeda kernel autonomous} is set to $(2,-0.3)$. Non-trivial equilibria $x_0=0.088$ (point $A$), $x_0=1.172$ (point $B$), $x_0=1.977$ (point $C$) are provided by the relation \eqref{Ikeda equilibrium equation}. Grey bands represent the regions where the sufficient stability conditions~\eqref{stability inequality} in the statement of the corollary are satisfied; the equilibria that correspond to the points $A$ and $C$ are hence stable.}
\label{equilibria for Ikeda} 
\end{figure}

Figure~\ref{equilibria for Ikeda} illustrates the statement of Corollary~\ref{suppl-Stability Ikeda continuous}.

\subsection{Fixed points of the reservoir map and their stability}

In this section we consider the discrete time TDR, we characterize its fixed points and establish sufficient conditions for their stability which, as we will show,  are analogous to the ones that we obtained for the continuous time case. 

We place ourselves in the autonomous regime, that is, $I(t) = 0$ in the time-delay differential equation \eqref{time delay equation} and $\mathbf{I} (t) = {\bf 0}_N $ in the associated recursion \eqref{vector discretized reservoir}. 
In this case, we can state the following proposition that shows that there is a bijective correspondence between the equilibria of the discrete and continuous time TDRs.
\begin{proposition}
\label{fixed point vs equilibrium} 
The point $x _0 \in \mathbb{R} $ is an equilibrium of the time-delay differential equation \eqref{time delay equation 2 recall}
in autonomous regime, that is when $ I(t)= 0$, if and only if
the vector $\mathbf{x} _0 :=  x_0 \mathbf{i} _N  $ is a fixed point of the $N$-dimensional discretized nonlinear time-delay reservoir
\begin{equation}
\label{}
\dot{\mathbf{x}} (t) = F( \mathbf{x} (t - 1), \mathbf{I} (t), \boldsymbol{\theta} )
\end{equation}
 in  autonomous regime, that is, when $\mathbf{I} (t) = {\bf 0} _N $.
\end{proposition}
{\bf{Proof.}}
Suppose  first that  $x _0 $ is an equilibrium of the time-delay differential equation \eqref{time delay equation 2 recall} and hence satisfies $x_0 = f( x_0 , 0 ,\boldsymbol{\theta} )$. In order to show that $F( \mathbf{x} (t - 1), \mathbf{0} _N , \boldsymbol{\theta} ) = \mathbf{x} _0 $, we evaluate the components of the right hand side of \eqref{discretized reservoir} at $\mathbf{x}_0 $ and obtain 
\begin{equation*}
(F(\mathbf{x} _0,{\bf 0} _N, \boldsymbol{ \theta }))_{i}={ e}^{-i \xi } x _0 +(1- { e}^{- \xi }) \sum^{i - 1}_{j = 0} { e}^{-j \xi } f( x _0 , 0,\boldsymbol{\theta} ) = { e}^{-i \xi } x _0 + (1 - { e}^{- \xi }) \sum^{i-1}_{j = 0} { e}^{-j \xi } x _0 = x _0,
\end{equation*}
as required. The proof of the converse implication is also straightforward using \eqref{discretized reservoir}.
$\square$

\medskip

We now provide  sufficient conditions for the stability of the fixed points  of the type $\mathbf{x}_0 = x_0 \mathbf{i} _N$ described in Proposition~\ref{fixed point vs equilibrium}. The asymptotic stability (respectively, stability) of those fixed points is guaranteed whenever the connectivity matrix $A(\mathbf{x} _0, \boldsymbol{ \theta })$ in \eqref{DxF} satisfies   $\rho (A(\mathbf{x} _0, \boldsymbol{ \theta })) < 1 $ (respectively, $ \rho (A(\mathbf{x} _0, \boldsymbol{ \theta })) \le 1)$). Since it is not possible to compute the eigenvalues of $A(\mathbf{x} _0, \boldsymbol{ \theta })$ for an arbitrary number of neurons $N$, we proceed by finding upper bounds of its spectral radius $\rho (A(\mathbf{x} _0, \boldsymbol{ \theta }))$. This can be done with the help of the Gershgorin disks theorem (see for instance Corollary 6.1.5 in~\cite{horn:matrix:analysis}) or by using a matrix norm $||| \cdot ||| $ and noting that $\rho(A(\mathbf{x} _0, \boldsymbol{ \theta })) \le  |||A(\mathbf{x} _0, \boldsymbol{ \theta })||| $. After a detailed study using all these possibilities we found that the best result is obtained by using  the maximum row sum matrix norm $|||A(\mathbf{x} _0, \boldsymbol{ \theta })|||_{\infty}$ defined in Section~\ref{Notation}, which allows us to formulate the following result.
\begin{theorem}
\label{stability discrete}
Let $\mathbf{x} _0 = x_0 \mathbf{i} _N $ be a fixed point of the $N$-dimensional recursion $\mathbf{x} (t) = F( \mathbf{x} (t-1), \mathbf{I} (t) , \boldsymbol{\theta} )$  in autonomous regime. Then, $\mathbf{x} _0 \in \mathbb{R}^N $ is asymptotically stable (respectively stable) if  
$|\partial _{{x}}f({ x} _0, { 0} , \boldsymbol{\theta}  )| < 1$ (respectively, $|\partial _{{x}}f({ x} _0, { 0} , \boldsymbol{\theta}  )| \le 1$).
\end{theorem}
{\bf Proof.} We first recall that the connectivity matrix \eqref{DxF} of the discretized  nonlinear TDR with $N$ virtual nodes is given by 
\begin{align}
\label{A}
&A(\mathbf{x} _0, \boldsymbol{ \theta }):=  
{\small\left(
\begin{array}{cccccc}
\Phi &    0 &\hdots &0&  { e^{- \xi }}\\
\Phi   { e^{- \xi }}&   \Phi &  \hdots &0& { e^{-2 \xi }}\\
 \Phi  e^{- 2\xi }&  \Phi  e^{- \xi }&  \hdots &0&  { e^{- 3\xi }}\\
  \vdots&\vdots&\vdots&&\vdots\\
  \Phi  e^{- (N-1)\xi }&  \Phi   e^{- (N-2)\xi }&   \hdots & \Phi e^{- \xi }& \Phi  + { e^{- N\xi }}
\end{array}
\right)},
\end{align}
where $\Phi : = (1-e^{- \xi })\cdot\partial _{{x}}f({ x} _0, { 0} , \boldsymbol{\theta}  )$  and recall that $\partial _{{x}}f({ x} _0, { 0} , \boldsymbol{\theta} ) $ is the first derivative of the nonlinear kernel $f$ in \eqref{time delay equation} with respect to the first argument and computed at the point $({\bf x} _0, {\bf 0}_N , \boldsymbol{\theta}  )$, with $ \xi  = \log (1 + d)$ and $d \in (0, 1]$  the Euler discretization step or, equivalently, the separation between the virtual neurons. We will use the notation $f'_{x_0}:=\partial _{{x}}f({ x} _0, {0} , \boldsymbol{\theta} ) $ in what follows.

We  proceed by finding sufficient conditions on $f'_{x_0}$ that guarantee that the spectral radius of $A(\mathbf{x} _0, \boldsymbol{ \theta })$ is bounded above by $1$. 
These conditions will be obtained by enforcing $|||A(\mathbf{x} _0, \boldsymbol{ \theta })|||_{\infty}  <1 $ and by recalling that 
\begin{equation}
\label{spectral radius}
\rho (A(\mathbf{x} _0, \boldsymbol{ \theta })) \le |||A(\mathbf{x} _0, \boldsymbol{ \theta })|||_{\infty} = \max _{1\le i \le N} \sum^{N}_{j = 1} |a _{ij} |.
\end{equation}
In the view of the rows of the matrix $A(\mathbf{x} _0, \boldsymbol{ \theta })$, it is clear that $ |||A(\mathbf{x} _0, \boldsymbol{ \theta })|||_{\infty}$ is given by the sum of the absolute values of one of the rows with numbers $1$, $N-1$, or $N$. We can hence write
\begin{equation}
\label{}
 |||A(\mathbf{x} _0, \boldsymbol{ \theta })|||_{\infty} = {\max} \Bigg\{ 
\begin{array}{cl} 
u:=&e^{- \xi } +  (1- e ^{- \xi })|f'_{x_0}|\\
 v:=&e^{- (N-1)\xi } +  (1- e ^{- (N-1)\xi })|f'_{x_0}|\\
  w:=&|e^{- N\xi } +  (1- e ^{- \xi })f'_{x_0}| + e^{- \xi }(1- e ^{- (N-1)\xi })|f'_{x_0}|
  \end{array}\Bigg\}.
\end{equation}
It can be easily verified that this expression can be split into two cases, namely
\begin{subnumcases}{|||A(\mathbf{x} _0, \boldsymbol{ \theta })|||_{\infty} = }
  \max \{ u, w\}, \label{case I} & if $|f'_{x_0} | \le 1$\\
    \max \{ v, w\}, \label{case II}  & if $|f'_{x_0} | \ge 1$
\end{subnumcases}
Additionally, by definition of the absolute value, we have two cases:
\begin{subnumcases}{|f'_{x_0}| = }
   f'_{x_0},& if $f'_{x_0}  \in [ 0, +\infty)$\label{case A} \\
    -f'_{x_0},  & if $f'_{x_0} \in (-\infty, 0)$ \label{case B} 
\end{subnumcases}
and hence
\begin{subnumcases}{|e^{- N\xi } +  (1- e ^{- \xi })f'_{x_0}| = }
   e^{- N\xi } +  (1- e ^{- \xi })f'_{x_0},& if $f'_{x_0}  \in \bigg[ -\dfrac{e^{- N\xi }}{1- e ^{- \xi }},+\infty\bigg)$ \label{case a} \\
    -e^{- N\xi } -  (1- e ^{- \xi })f'_{x_0},  & if $f'_{x_0}  \in \bigg(-\infty, -\dfrac{e^{- N\xi }}{1- e ^{- \xi }}\bigg)$. \label{case b} 
\end{subnumcases}
We now consider in detail all the possible combinations of cases that provide the conditions on $f'_{x_0}$ that ensure stability by enforcing that $|||A(\mathbf{x} _0, \boldsymbol{ \theta })|||_{\infty}  <1 $.

\noindent{\bf Case I \eqref{case I}, \eqref{case A}, \eqref{case a}.} 

\noindent On one hand, \eqref{case I}, \eqref{case A}, \eqref{case a} give that $f'_{x_0} \in [0, 1)$. On the other hand  \eqref{case I} amounts to
\begin{subnumcases}{|||A(\mathbf{x} _0, \boldsymbol{ \theta })|||_{\infty} = \max \{ u, w\} = }
   u,& if $f'_{x_0}  \in(-\infty,  1]$\\
    w, & if $f'_{x_0}  \in ( 1, +\infty).$
\end{subnumcases}
Notice that since $f'_{x_0} \in [0, 1)$, then these cases reduce to:
\begin{equation}
\label{}
|||A(\mathbf{x} _0, \boldsymbol{ \theta })|||_{\infty} = u, 
\end{equation}
and hence
\begin{equation*}
\label{}
\rho (A(\mathbf{x} _0, \boldsymbol{ \theta }))\le u = e^{- \xi } +(1- e^{- \xi }) f'_{x_0}<1 \enspace \Longrightarrow f'_{x_0} \in (-\infty, 1).
\end{equation*}
We  finally write that in this case stability is guaranteed whenever
\begin{equation}
\label{I1}
f'_{x_0} \in I1:= [0, 1).
\end{equation}
\noindent{\bf Case II \eqref{case I}, \eqref{case A}, \eqref{case b}.} 

\noindent On one hand, \eqref{case I}, \eqref{case A} imply that $f'_{x_0} \in [0, 1)$ but by \eqref{case b} it is required at the same time  that $f'_{x_0}  \in (-\infty, -\dfrac{e^{- N\xi }}{1- e ^{- \xi }})$ which immediately yields in this case:
\begin{equation}
\label{I2}
f'_{x_0} \in I2:=  \emptyset. 
\end{equation}
\noindent{\bf Case III \eqref{case I}, \eqref{case B}, \eqref{case a}.} 

\noindent On one hand, \eqref{case I}, \eqref{case B}, \eqref{case a} imply that \begin{equation}
\label{Case III initial cond}
f'_{x_0} \in \Big[\max\Big\{-1, -\dfrac{e^{- N\xi }}{1- e ^{- \xi }}\Big\}, 0\Big).
\end{equation}
On the other hand, the condition that defines \eqref{case I} amounts to
\begin{subnumcases}{\!\!\!\!\!\!\!\!\!\!\!\!\!\!\!\!\!\!\!\!\!\!\!\!\!\!\!\!|||A(\mathbf{x} _0, \boldsymbol{ \theta })|||_{\infty} = \max \{ u, w\} =}
   u,& if $f'_{x_0} \in \bigg[-\dfrac{e^{- \xi } - e^{-N \xi }}{|2-3 e^{- \xi }+ e^{-N \xi }|},   \dfrac{e^{- \xi } - e^{-N \xi }}{|2-3 e^{- \xi }+ e^{-N \xi }|}\bigg]$\label{case III 1} \\
    w, & $\begin{array}{l}
{\rm if } \enspace f'_{x_0} \in  \bigg(-\infty, -\dfrac{e^{- \xi } - e^{-N \xi }}{|2-3 e^{- \xi }+ e^{-N \xi }|}\bigg) \\ \cup    \bigg(\dfrac{e^{- \xi } - e^{-N \xi }}{|2-3 e^{- \xi }+ e^{-N \xi }|}, +\infty\bigg)
\end{array}. $ \label{case III 2}
\end{subnumcases}
Notice now that at the same time we require that $|||A(\mathbf{x} _0, \boldsymbol{ \theta })|||_{\infty}<1$. Hence for the case \eqref{case III 1} we have
\begin{equation}
\label{rho A case III u}
 \rho (A(\mathbf{x} _0, \boldsymbol{ \theta })) \le |||A(\mathbf{x} _0, \boldsymbol{ \theta })|||_{\infty} = u=e^{- \xi } -  (1- e ^{- \xi })f'_{x_0}< 1\Longrightarrow f'_{x_0}\in (-1,+\infty)
\end{equation}
which put together with the conditions for $f'_{x_0}$ in \eqref{case III 1}, \eqref{rho A case III u}, and  \eqref{Case III initial cond} yields
\begin{equation}
\label{I3}
f'_{x_0} \in I3:=  \Bigg(\max \Bigg\{-1, -\dfrac{e^{- N\xi }}{1- e ^{- \xi }}, -\dfrac{e^{- \xi } - e^{-N \xi }}{|2-3 e^{- \xi }+ e^{-N \xi }|}\Bigg\}, 0\Bigg). 
\end{equation}
We now consider the case \eqref{case III 2}  and in an analogous way using \eqref{spectral radius} we have that
\begin{align}
\label{rho A case III w}
 \rho (A(\mathbf{x} _0, \boldsymbol{ \theta })) \le |||A(\mathbf{x} _0, \boldsymbol{ \theta })|||_{\infty} = &w=e^{- N\xi } +  (1- 2 e ^{- \xi } + e ^{- N\xi })f'_{x_0}< 1 \Longrightarrow \nonumber \\&f'_{x_0} \in \bigg(-\dfrac{1-e^{-N \xi }}{|1-2e^{- \xi }+ e ^{- N\xi }|},\dfrac{1-e^{-N \xi }}{|1-2e^{- \xi }+ e ^{- N\xi }|}\bigg)
\end{align}
If we put together the conditions for $f'_{x_0}$ in \eqref{case III 2}, \eqref{rho A case III w}, and  \eqref{Case III initial cond} we obtain
\begin{equation}
\label{I4}
f'_{x_0} \in I4:=  \Bigg(\max \Bigg\{-1, -\dfrac{e^{- N\xi }}{1- e ^{- \xi }},-\dfrac{1-e^{-N \xi }}{|1-2e^{- \xi }+ e ^{- N\xi }|}\Bigg\}, -\dfrac{e^{- \xi } - e^{-N \xi }}{|2-3 e^{- \xi }+ e^{-N \xi }|}\Bigg). 
\end{equation}
\noindent{\bf Case IV \eqref{case I}, \eqref{case B}, \eqref{case b}.} 

\noindent On  one hand, \eqref{case I}, \eqref{case B}, \eqref{case b} imply that 
\begin{equation}
\label{Case IV initial cond*}
f'_{x_0} \in \Big(-1, -\dfrac{e^{- N\xi }}{1- e ^{- \xi }}\Big).
\end{equation}
On the other hand, the case \eqref{case I} amounts to
\begin{subnumcases}{|||A(\mathbf{x} _0, \boldsymbol{ \theta })|||_{\infty} = \max \{ u, w\} = }
   u,& if $f'_{x_0}  \in \bigg[   -\dfrac{1 + e^{-(N-1) \xi }}{1 - e^{-(N-1) \xi }},+\infty\bigg)$\label{case IV 1}\\
   w,& if $f'_{x_0}  \in \bigg(-\infty,   -\dfrac{1 + e^{-(N-1) \xi }}{1 - e^{-(N-1) \xi }}\bigg)$\label{case IV 2}.
\end{subnumcases}
It can be easily verified that the condition defining \eqref{case IV 2} is incompatible with \eqref{Case IV initial cond*} since the relation $-1\ge-\dfrac{1 + e^{-(N-1) \xi }}{1 - e^{-(N-1) \xi }}$ holds for any $\xi \in (0, 1]$ and $N \in \mathbb{N} $. We hence conclude  using the case \eqref{case IV 1}  that
\begin{equation*}
\label{rho A case IV u}
 \rho (A(\mathbf{x} _0, \boldsymbol{ \theta })) \le |||A(\mathbf{x} _0, \boldsymbol{ \theta })|||_{\infty} = u=e^{- \xi } -  (1- e ^{- \xi } )f'_{x_0}< 1 \Longrightarrow f'_{x_0} \in (-1, +\infty).
\end{equation*}
Hence, together with \eqref{Case IV initial cond*} this amounts to
\begin{equation}
\label{Case IV initial cond}
f'_{x_0} \in I5:=\Big(-1, -\dfrac{e^{- N\xi }}{1- e ^{- \xi }}\Big).
\end{equation}
\noindent{\bf Case V \eqref{case II}, \eqref{case A}, \eqref{case a}.} 

\noindent On one hand, \eqref{case II}, \eqref{case A}, \eqref{case a} imply that 
\begin{equation}
\label{Case V initial cond}
f'_{x_0} \in [1,+\infty].
\end{equation}
On the other hand, the case \eqref{case I} amounts to
\begin{subnumcases}{|||A(\mathbf{x} _0, \boldsymbol{ \theta })|||_{\infty} = \max \{ v, w\} = }
   v,& if $f'_{x_0}  \in (-\infty, -1]$  \label{case V 1} \\
   w,&  if $f'_{x_0}  \in  (-1,+\infty)$\label{case V 2}.
\end{subnumcases}
Due to \eqref{Case V initial cond}
\begin{equation*}
\label{rho A case V w}
 \rho (A(\mathbf{x} _0, \boldsymbol{ \theta })) \le |||A(\mathbf{x} _0, \boldsymbol{ \theta })|||_{\infty} = w=e^{- N\xi } +  (1- e ^{- N\xi } )f'_{x_0}< 1 \Longrightarrow f'_{x_0} \in (-\infty, 1),
\end{equation*}
which yields 
\begin{equation}
\label{I6}
f'_{x_0} \in I6:=  \emptyset. 
\end{equation}
\noindent{\bf Case VI \eqref{case II}, \eqref{case A}, \eqref{case b}.} 

\noindent The conditions \eqref{case II}, \eqref{case A}, and  \eqref{case b} immediately yield  
\begin{equation}
\label{I7}
f'_{x_0} \in I7:=  \emptyset. 
\end{equation}
\noindent{\bf Case VII \eqref{case II}, \eqref{case B}, \eqref{case a}.} 

\noindent On  one hand, \eqref{case II}, \eqref{case B}, \eqref{case a} imply that 
\begin{equation}
\label{Case VII initial cond}
f'_{x_0} \in \Big[-\dfrac{e^{- N\xi }}{1- e ^{- \xi }},-1\Big].
\end{equation}
On the other hand, the case \eqref{case II} amounts to
\begin{subnumcases}{|||A(\mathbf{x} _0, \boldsymbol{ \theta })|||_{\infty} = \max \{ v, w\} = }
   v,& if $f'_{x_0}  \in \bigg(-\infty,  -\dfrac{ e^{-(N-1) \xi }}{2-e^{-(N-1) \xi }} \bigg]$\label{case VII 1}, \\
   w,&  if $f'_{x_0}  \in \bigg(  -\dfrac{ e^{-(N-1) \xi }}{2-e^{-(N-1) \xi }} ,+\infty\bigg).$\label{case VII 2} 
\end{subnumcases}
Hence, for the case \eqref{case VII 1} we have
\begin{equation*}
\label{rho A case V|| v}
 \rho (A(\mathbf{x} _0, \boldsymbol{ \theta })) \le |||A(\mathbf{x} _0, \boldsymbol{ \theta })|||_{\infty} = v=e^{- (N-1)\xi } -  (1- e ^{-(N-1) \xi })f'_{x_0}< 1\Longrightarrow f'_{x_0}\in (-1,+\infty)
\end{equation*}
which due to \eqref{Case VII initial cond} gives
\begin{equation}
\label{I8}
f'_{x_0} \in I8:=  \emptyset. 
\end{equation}
We now consider the case \eqref{case VII 2}  and in an analogous way we have that
\begin{align}
\label{rho A case V|| w}
 \rho (A(\mathbf{x} _0, \boldsymbol{ \theta })) \le |||A(\mathbf{x} _0, \boldsymbol{ \theta })|||_{\infty} = &w=e^{- N\xi } +  (1- 2 e ^{- \xi } + e ^{- N\xi })f'_{x_0}< 1 \Longrightarrow \nonumber \\&f'_{x_0} \in \bigg(-\dfrac{1-e^{-N \xi }}{|1-2e^{- \xi }+ e ^{- N\xi }|},\dfrac{1-e^{-N \xi }}{|1-2e^{- \xi }+ e ^{- N\xi }|}\bigg).
\end{align}
Hence by \eqref{rho A case V|| w}, \eqref{case VII 2} and \eqref{Case VII initial cond} we can write\begin{equation}
\label{I9}
f'_{x_0} \in I9:=  \bigg( \max \bigg\{-\dfrac{e^{- N\xi }}{1- e ^{- \xi }}, -\dfrac{ e^{-(N-1) \xi }}{2-e^{-(N-1) \xi }}, -\dfrac{1-e^{-N \xi }}{|1-2e^{- \xi }+ e ^{- N\xi }|}\bigg\}, -1 \bigg] = \emptyset,
\end{equation}
where the last equality follows from the fact that $-\dfrac{ e^{-(N-1) \xi }}{2-e^{-(N-1) \xi }} \ge -1$, for any $\xi \in (0,1]$ and $N \in \mathbb{N} $.

\noindent{\bf Case VIII \eqref{case II}, \eqref{case B}, \eqref{case b}.} 

\noindent On one hand, \eqref{case II}, \eqref{case B}, \eqref{case b} imply that 
\begin{equation}
\label{Case VIII initial cond}
f'_{x_0} \in \Big(-\infty, \min \bigg\{-\dfrac{e^{- N\xi }}{1- e ^{- \xi }},-1\bigg\}\Big).
\end{equation}
On the other hand, the case \eqref{case II} amounts to
\begin{subnumcases}{|||A|||_{\infty} = \max \{ v, w\} = }
   v,& if $f'_{x_0}  \in \bigg[-\dfrac{ 1 + e^{- \xi }}{1 - e^{- \xi }},+\infty \bigg)$\label{case VIII 1} \\
   w,&  if $f'_{x_0}  \in \bigg(  -\infty,-\dfrac{ 1 + e^{- \xi }}{1 - e^{- \xi }}\bigg).$\label{case VIII 2} 
\end{subnumcases}
Hence for the case \eqref{case VIII 1} we have
\begin{equation*}
\label{rho A case V|I| v}
 \rho (A(\mathbf{x} _0, \boldsymbol{ \theta })) \le |||A(\mathbf{x} _0, \boldsymbol{ \theta })|||_{\infty} = v=e^{- (N-1)\xi } -  (1- e ^{-(N-1) \xi })f'_{x_0}< 1\Longrightarrow f'_{x_0}\in (-1,+\infty)
\end{equation*}
which due to \eqref{Case VIII initial cond} implies that
\begin{equation}
\label{I10}
f'_{x_0} \in I10:=  \emptyset. 
\end{equation}
We now consider the case \eqref{case VIII 2}  and in an analogous way  we have that
\begin{equation}
\label{rho A case V|I| w}
 \rho (A(\mathbf{x} _0, \boldsymbol{ \theta })) \le |||A(\mathbf{x} _0, \boldsymbol{ \theta })|||_{\infty} = w=-e^{- N\xi } +  (-1+ e ^{- N\xi })f'_{x_0}< 1 \Longrightarrow f'_{x_0} \in \bigg(-\dfrac{1+e^{-N \xi }}{1-e^{-N \xi }},+\infty \bigg).
\end{equation}
Hence by \eqref{rho A case V|I| w}, \eqref{case VIII 2} and \eqref{Case VIII initial cond} we can write\begin{equation}
\label{I11}
f'_{x_0} \in I11:=  \bigg(-\dfrac{1+e^{-N \xi }}{1-e^{-N \xi }},  \min \bigg\{-\dfrac{e^{- N\xi }}{1- e ^{- \xi }}, -1, -\dfrac{ 1 + e^{- \xi }}{1 - e^{- \xi }} \bigg\}\bigg) = \emptyset,
\end{equation}
for any $\xi \in (0,1]$ and $N \in \mathbb{N} $.

\medskip

Finally, we put together all the non-empty intervals provided by the  Cases I-VIII that guarantee that when $f'_{x_0} $ belongs to them, then the fixed point $\mathbf{x} _0 $ is stable:
\begin{align}
\label{non-empty intervals for all cases}
f'_{x_0}\in \Big(-1, -\dfrac{e^{- N\xi }}{1- e ^{- \xi }}\Big) &\cup \Bigg(\max \Bigg\{-1, -\dfrac{e^{- N\xi }}{1- e ^{- \xi }},-\dfrac{1-e^{-N \xi }}{|1-2e^{- \xi }+ e ^{- N\xi }|}\Bigg\}, -\dfrac{e^{- \xi } - e^{-N \xi }}{|2-3 e^{- \xi }+ e^{-N \xi }|}\Bigg) \nonumber \\
&\cup  \Bigg(\max \Bigg\{-1, -\dfrac{e^{- N\xi }}{1- e ^{- \xi }}, -\dfrac{e^{- \xi } - e^{-N \xi }}{|2-3 e^{- \xi }+ e^{-N \xi }|}\Bigg\}, 0\Bigg) \cup   [0,1).
\end{align}
Now, it is easy to see that 
\begin{align*}
\label{}
&\max \Bigg\{-1, -\dfrac{e^{- N\xi }}{1- e ^{- \xi }}, -\dfrac{e^{- \xi } - e^{-N \xi }}{|2-3 e^{- \xi }+ e^{-N \xi }|}\Bigg\} \ge -1, \nonumber \\
&\max \Bigg\{-1, -\dfrac{e^{- N\xi }}{1- e ^{- \xi }},-\dfrac{1-e^{-N \xi }}{|1-2e^{- \xi }+ e ^{- N\xi }|}\Bigg\} \ge -1,
\end{align*}
and hence the condition \eqref{non-empty intervals for all cases} reduces to
\begin{equation*}
\label{}
|f'_{x_0}| < 1.
\end{equation*}
Notice that the proof remains valid for the case of (not necessary asymptotic) stability in the statement of the theorem. In this case we require
\begin{equation*}
\label{spectral radius_equality}
\rho (A(\mathbf{x} _0, \boldsymbol{ \theta })) \le |||A(\mathbf{x} _0, \boldsymbol{ \theta })|||_{\infty} = \max _{1\le i \le N} \sum^{N}_{j = 1} |a _{ij} | \le1,
\end{equation*}
which results in  the condition
\begin{equation*}
\label{}
|f'_{x_0}| \le 1,
\end{equation*}
as required.
$\square$

\medskip

The following theorem provides the characteristic polynomial of the connectivity matrix $A(\mathbf{x} _0, \boldsymbol{ \theta })$ and an explicit expression for its spectral radius under some conditions. Another way to find upper bounds for this spectral radius would consist of using the Cauchy bound~\cite{Rahman:Schmeisser} of this polynomial. In our experience this approach produces mediocre results in comparison with the statement in Theorem~\ref{stability discrete}.

\begin{theorem}
Let $A(\mathbf{x} _0, \boldsymbol{ \theta })$ be the connectivity matrix of the reservoir map in \eqref{DxF}. Define $\Phi :=(1- e^{- \xi })f'_{x_0}$ and let $\left\{\lambda_1, \dots, \lambda_N  \right\} $ be the roots of the polynomial  equation
\begin{equation}
\label{polynomial  equation}
\lambda ^N -\left( \dfrac{e^{-N \xi }}{\Phi}+N\right) \lambda ^{N-1} + \sum^{N-2}_{j = 0}\left( \begin{array}{c}N \\ j\end{array}\right) (-1)^{N-j} \lambda ^j =0.
\end{equation}
Then $\rho (A(\mathbf{x} _0, \boldsymbol{ \theta })) = \max\left\{ |\Phi \lambda _1 |, \dots,   |\Phi \lambda _N |\right\} $. Moreover, if \eqref{polynomial  equation} has a root $\lambda$ such that $\lambda>1$, then $\rho (A(\mathbf{x} _0, \boldsymbol{ \theta })) =\lambda |\Phi |$ necessarily.  
\end{theorem}
\noindent{\bf Proof.} Define first $B:=\dfrac{1}{\Phi}A(\mathbf{x} _0, \boldsymbol{ \theta })$. We then write
\begin{equation*}
\label{}
B:= {\small\left(
\begin{array}{cccccc}
1 &    0 &\hdots &0& \dfrac{ e^{- \xi }}{ \Phi }\\
   { e^{- \xi }}&   1 &  \hdots &0& \dfrac{ e^{-2 \xi }}{ \Phi }\\
   e^{- 2\xi }&    e^{- \xi }&  \hdots &0&  \dfrac{ e^{- 3\xi }}{ \Phi }\\
  \vdots&\vdots&\vdots& &\vdots\\
    e^{- (N-1)\xi }&     e^{- (N-2)\xi }&   \hdots &  e^{- \xi }& 1  +\dfrac{ e^{- N\xi }}{\Phi }
\end{array}
\right)}.
\end{equation*}
Let $\mathbf{v} = (v _1 , \dots, v _N )$ be an eigenvector of $B$ with eigenvalue $\lambda$, that is, $B\mathbf{v} = \lambda\mathbf{v} $. This equality can be rewritten as
\begin{equation}
\label{system for eigenvalues of B}
\left\{
\begin{array}{lcl}
v _1 &=& \dfrac{e^{- \xi }}{(\lambda - 1)\Phi } v _N,  \\
v _2 &=& \dfrac{e^{- 2\xi }}{(\lambda - 1)\Phi } v _N \dfrac{\lambda }{ \lambda -1 },\\
v _3 &=& \dfrac{e^{- 3\xi }}{(\lambda - 1)\Phi} v _N \dfrac{\lambda ^2 }{( \lambda -1 ) ^2 },\\
\vdots&\vdots&\ \ \ \ \ \ \ \ \ \vdots\\
v _{ N - 1} &=& \dfrac{e^{- (N-1)\xi }}{(\lambda - 1)\Phi } v _N \dfrac{\lambda ^{N-2} }{ (\lambda -1) ^{N-2}},
\end{array} 
\right.
\end{equation}
together with the identity
\begin{equation*}
\label{}
\sum^{N}_{j = 1} b _{Nj} v _j = \lambda v_N,  
\end{equation*}
which is equivalent to 
\begin{equation*}
\label{}
\dfrac{v_N}{(\lambda - 1)\Phi }e^{-N \xi }\left(1+ \dfrac{ \lambda }{\lambda - 1 } +  \dfrac{ \lambda ^2 }{(\lambda -1) ^2 } +\cdots + \dfrac{ \lambda^{N-2} }{(\lambda -1)^{N-2}} \right)  + v _N \left(1 + \dfrac{e^{-N \xi }}{\Phi}\right) = \lambda v _N,
\end{equation*}
or, equivalently, to
\begin{equation*}
\label{}
v _N  \left(\dfrac{e^{-N \xi }}{ \Phi } \dfrac{\lambda^{N-1} - ( \lambda -1 )^{N-1}}{( \lambda - 1)^{N-1}}+ ( \lambda - 1 ) + \dfrac{e^{-N \xi }}{ \Phi }\right) = 0. 
\end{equation*}
If we assume that $b _N \neq 0$ and $ \lambda \neq 1$, this amounts to
\begin{equation*}
\label{}
\dfrac{( \lambda -1 )^{N}}{\lambda^{N-1}} = \dfrac{e^{-N \xi }}{ \Phi }.
\end{equation*}
Consequently, the eigenvalues of $A$ are given by the roots $ \left\{ \lambda _1 , \dots, \lambda _N \right\} $ of the polynomial equation
\begin{equation*}
\label{}
\lambda ^N -\left( \dfrac{e^{-N \xi }}{\Phi}+N\right) \lambda ^{N-1} + \sum^{N-2}_{j = 0}\left( \begin{array}{c}N \\ j\end{array}\right) (-1)^{N-j} \lambda ^j =0
\end{equation*}
and hence \begin{equation*}
\label{}
\rho (A(\mathbf{x} _0, \boldsymbol{ \theta })) = \max\left\{ |\Phi \lambda _1 |, \dots,   |\Phi \lambda _N |\right\}.
\end{equation*}
If there exists some $\lambda>1$, the expressions \eqref{system for eigenvalues of B} show that the eigenvector $\mathbf{v}$ of $B$ (or of $-B$ if $\Phi< 0$) can be chosen positive. In that situation Corollary 8.1.30 in \cite{horn:matrix:analysis} guarantees that $\rho(B) =\lambda  $ with $\lambda>1$  the eigenvalue corresponding to $\mathbf{v} $ and hence $\rho (A(\mathbf{x} _0, \boldsymbol{ \theta })) = | \Phi \lambda |$ as required. $\square$

\section{Robustness of the empirical tests with respect to the choice of nonlinear kernel}

In this section we show the robustness of the empirical results in Sections~\ref{Optimal performance: stability and unimodality}  and~\ref{Discussion} of the paper with respect to the choice of the nonlinear kernel used in the construction of the TDR. 

\begin{figure}[ht]
\centering
\begin{tikzpicture}
    \node[anchor=south west,inner sep=0] at (0,0) {\hspace{-2.5cm}\includegraphics[scale=.3]{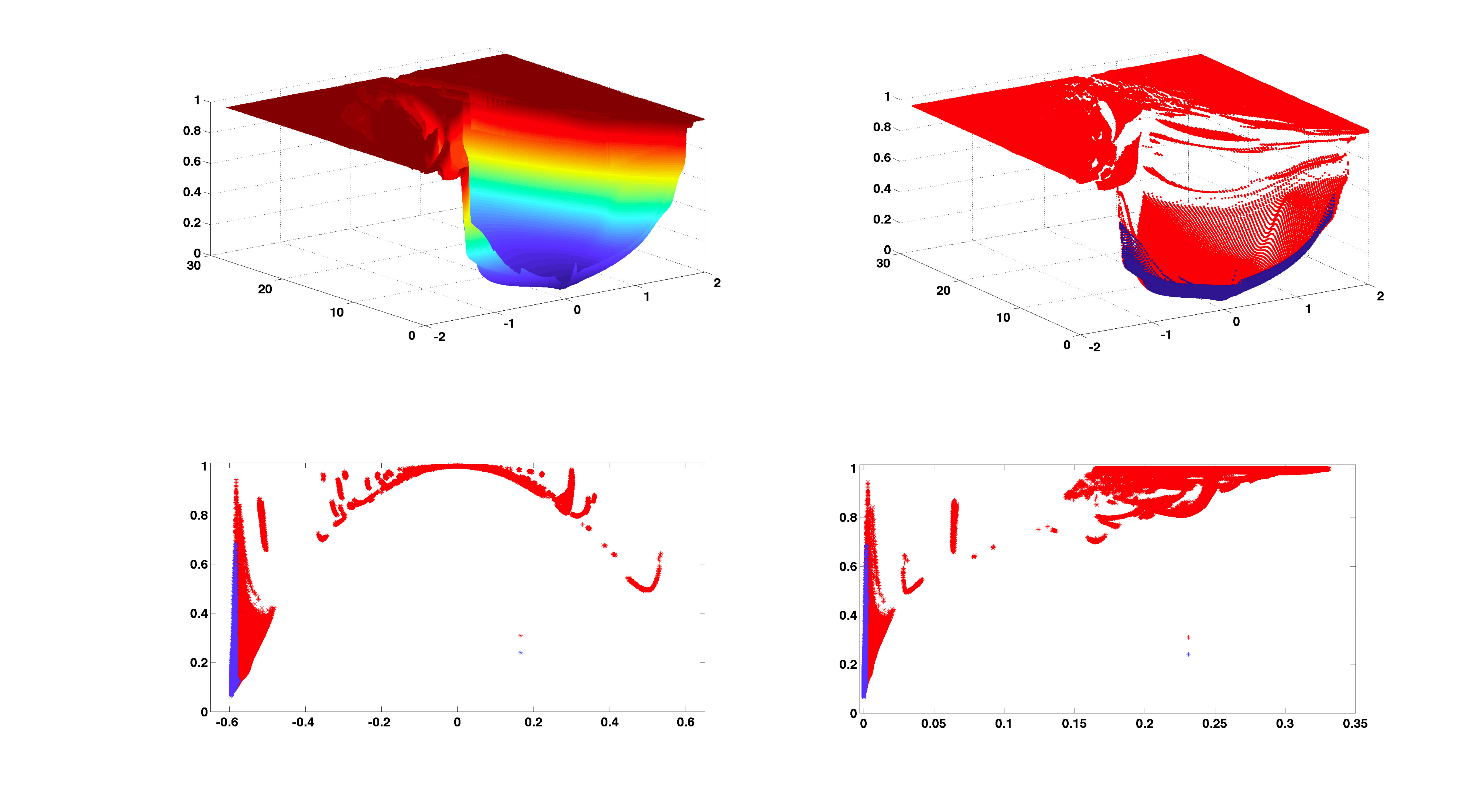}};
          \node at (5.4,2.43) {\tiny Bimodality};          
          \node at (5.46,2.2) {\tiny Unimodality};           
          \node at (14.58,2.43) {\tiny Bimodality};          
          \node at (14.64,2.2) {\tiny Unimodality}; 
           \node at (12.64,0.9) {\scriptsize Variance of the reservoir output}; 
          \node at (12.64,5.4) {\scriptsize Influence of the reservoir output variance on the normalized}; 
          \node at (12.64,5.1) {\scriptsize mean square error (NMSE) in the $3$-lag quadratic memory task};  \node[label=below:\rotatebox{90}{\scriptsize NMSE}] at (8.8,3.7) {};
             \node[label=below:\rotatebox{90}{\scriptsize NMSE}] at (-0.16,3.7) {}; 
             \node at (4,0.9) {\scriptsize Mean of the reservoir output}; 
          \node at (4,5.4) {\scriptsize Influence of the reservoir output mean on the normalized}; 
          \node at (4,5.1) {\scriptsize mean square error (NMSE) in the $3$-lag quadratic memory task};
           \node[label=below:\rotatebox{90}{\scriptsize NMSE}] at (9.2,9.4) {};
           \node[label=below:\rotatebox{90}{\scriptsize NMSE}] at (-0.19,9.4) {};
           \node at (12.64,11.3) {\scriptsize Influence of the input mask on on the normalized}; 
          \node at (12.64,11) {\scriptsize mean square error (NMSE) in the $3$-lag quadratic memory task};
          \node at (4,11.3) {\scriptsize Influence of the input mask on on the normalized}; 
          \node at (4,11) {\scriptsize mean square error (NMSE) in the $3$-lag quadratic memory task};
          \node[label=below:\rotatebox{11}{\scriptsize Mean of the input mask}] at (5.7,7.25) {};
          \node[label=below:\rotatebox{-17}{\scriptsize Variance of the input mask}] at (1.7,7.63) {};
          \node[label=below:\rotatebox{10}{\scriptsize Mean of the input mask}] at (14.6,7.10) {};
          \node[label=below:\rotatebox{-23}{\scriptsize Variance of the input mask}] at (10.5,7.9) {};
\end{tikzpicture}
\caption{Behavior of the performance of a Mackey-Glass based reservoir in a quadratic memory task as a function of the mean and the variance of the input mask. The modification of any of these two parameters influences how the reservoir dynamics separates from the stable equilibrium. The top panels show how the performance degrades very quickly as soon as the  mean and the variance  of the input mask (and hence of the input signal) separate from zero. The bottom panels depict the reservoir performance as a function of the various output means and variances obtained when changing the input means and variances. In the top right and bottom panels we have indicated with  red markers the cases in which the reservoir visits the stability basin of a contiguous equilibrium hence showing how unimodality is associated to optimal performance. }
\label{unimodality_via_simulatons_MG}
\end{figure}
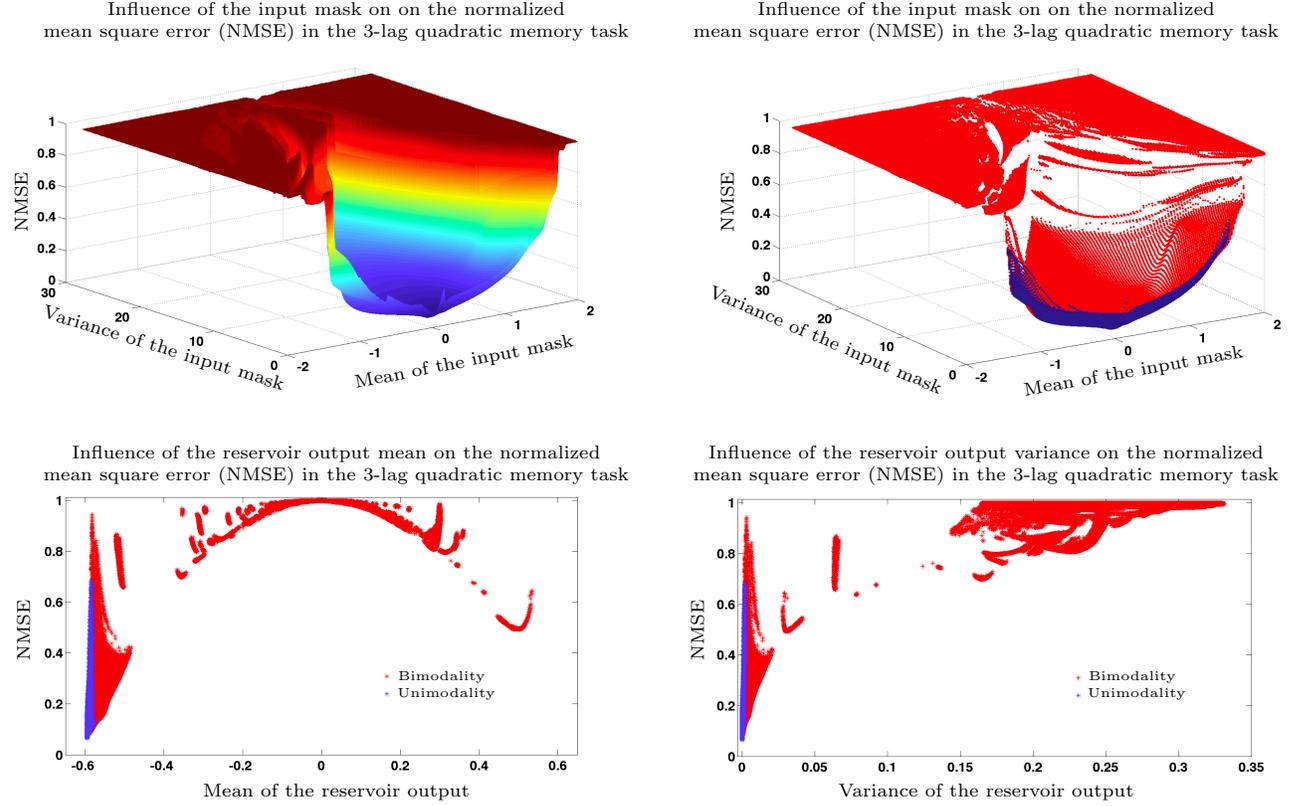

First, in Section~\ref{Optimal performance: stability and unimodality} we carried out an experiment using the Ikeda kernel and a quadratic memory task that showed that optimal performance is obtained when the input mean and variance are tuned so that the dynamics of the reservoir takes place in the neighborhood of a stable steady state and making sure that multimodality is  avoided. We have repeated here the same experiment but, this time, using the Mackey-Glass kernel. More specifically, we nonsider a TDR  with $N=20$ neurons $d=0.943$, $\gamma=4.7901 $, $\eta= 1.3541 $, and $p=2 $. As we explained in Corollary~\ref{Stability MG continuous}, with these parameter values  $x _0 = -\sqrt{\eta - 1} = -0.5951$ is an equilibrium that satisfies the sufficient conditions for asymptotic stability. 
In order to verify that the optimal performance is obtained  when the RC operates in a neighborhood of that stable equilibrium, we study the normalized mean square  error (NMSE) exhibited by a TDR initialized at $x _0=-0.5951 $ when we present to it a quadratic memory task. More specifically, we inject in a TDR under study an independent and identically normally distributed signal $z (t)$ with mean zero and variance $10^{-4} $ and we then train a linear readout $W_{{\rm  out}}$ (obtained with a ridge penalization of $\lambda=10^{-15}$) in order to recover the quadratic function $z (t-1)^2+  z (t-2)^2+z (t-3)^2$ out of the reservoir output.  The top left panel in Figure~\ref{unimodality_via_simulatons_MG} shows how the NMSE behaves as a function of the mean and the variance of the input mask $ {\bf c} $. It is clear that by modifying any of these two parameters we control how far the reservoir dynamics separates from the stable equilibrium, which we quantitatively evaluate in the two bottom panels by representing the RC performance in terms of the mean and the variance of the resulting reservoir output. Both panels depict how the injection of a signal slightly shifted in mean or with a sufficiently high variance results in reservoir outputs that separate from the stable equilibrium and in a severely degraded performance. An important factor in this deterioration seems to be the multi modality, that is, if the shifting in mean or the input signal variance are large enough then the reservoir output visits the stability basin of the other stable point placed at $x _0=\sqrt{\eta - 1} = 0.5951 $; in the top right and bottom panels we have marked with red color the values for which bimodality has occurred so that the negative effect of this phenomenon is noticeable. 
%

\begin{figure}[ht]
\hspace*{-2cm}\includegraphics[scale=.3]{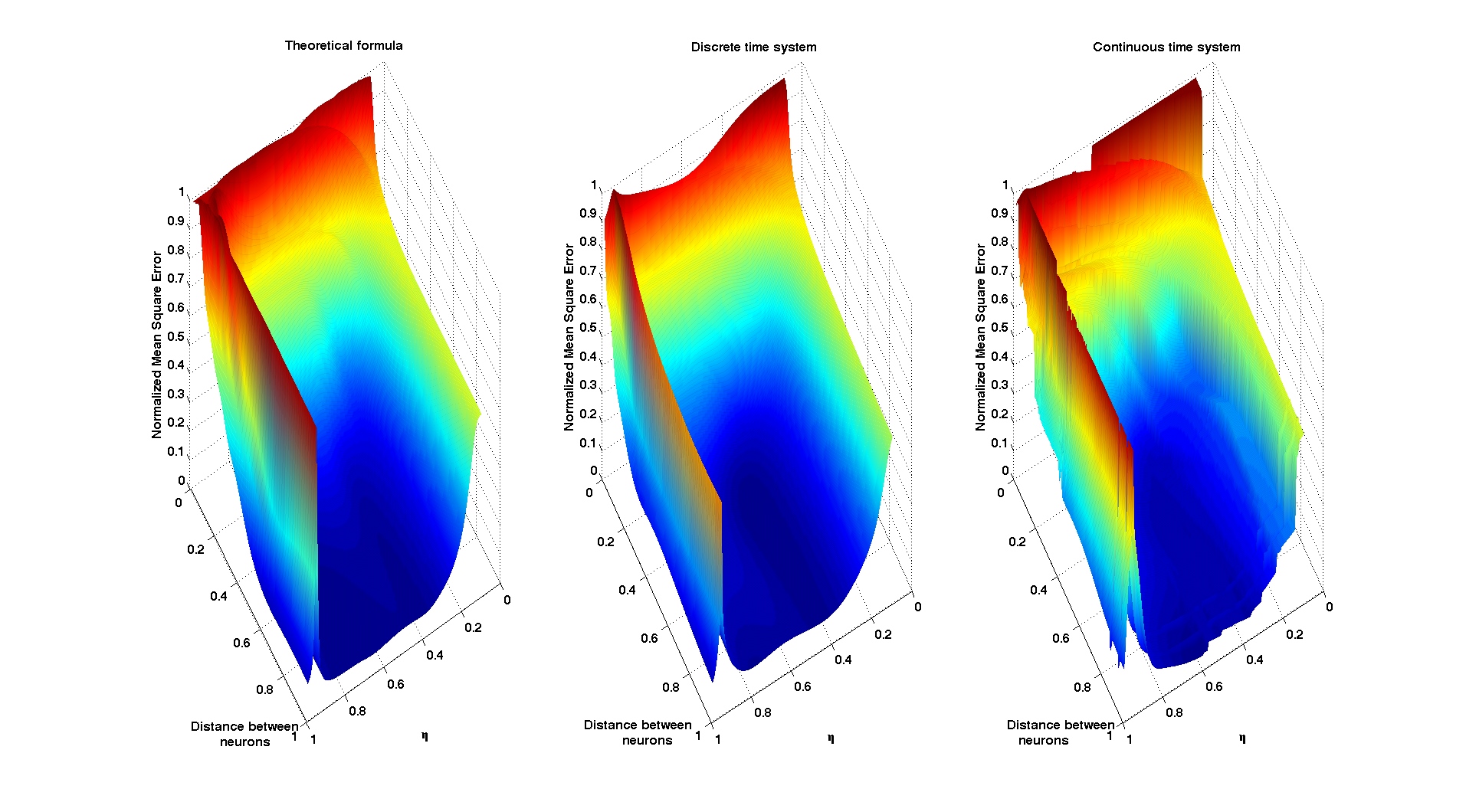}
\caption{Error surfaces exhibited by an Ikeda  based reservoir computer in a 3-lag quadratic memory task, as a function of the distance between neurons and the parameter $\eta$.  The points in the surfaces of the middle and right panels are the result of Monte Carlo evaluations of the NMSE exhibited by the discrete and continuous time TDRs, respectively.  The left panel was constructed using the formula~\eqref{capacity formula} that is obtained as a result of modeling the reservoir with an approximating VAR(1) process.}
\label{Performace_Ikeda}
\end{figure}

Second, in Section~\ref{Discussion} we used a Mackey-Glass based reservoir to compare the empirical performance surfaces in terms of various parameters with that coming from the formula~\eqref{capacity formula} that was obtained as a result of modeling the reservoir with an approximating VAR(1) process. In this section we have repeated the same exercise with an Ikeda based reservoir in order to show that the formula~\eqref{capacity formula} produces in this case results of comparable quality. The outcome of this experiment are contained in Figure~\ref{Performace_Ikeda} where we represent the normalized mean square error as a function of the distance between neurons and the feedback gain $\eta $. The other fixed parameter values used are $\gamma= 0.523 $ and $\phi=0.3106 $; the reservoir was constructed using 20 neurons and we presented to it the three-lag quadratic memory task corresponding to the diagonal matrix $Q$ with diagonal entries given by the vector $(0,1,1,1)$. The optimal output mask ${\mathbf W}_{{\rm out}}$ was computed using a ridge regression with $\lambda=10^{-15}$. As it can be seen in the figure, we restrict the values of the parameter $\eta$ to the interval $[ 0 ,1] $ which ensures, using one of the results in Corollary~\ref{suppl-Stability Ikeda continuous}, that the TDR exhibits for each value of $\eta $ a unique equilibrium (unimodality is hence guaranteed) that is always stable. The TDR is always initialized at that stable configuration.
\end{appendices}
\addcontentsline{toc}{section}{Bibliography}
\bibliographystyle{wmaainf}
\bibliography{/Users/Lyudmila/Mendeley/GOLibrary}
\end{document}